\documentclass[accepted]{uai2025} 
                        
\usepackage[american]{babel}

\usepackage{natbib} 
\usepackage{mathtools} 
\usepackage{booktabs} 
\usepackage{tikz} 
\usepackage{textcase}

\usepackage{graphicx}
\usepackage{apptools}
\usepackage[flushleft]{threeparttable}
\usepackage{array,booktabs,makecell}
\usepackage{multirow}
\usepackage{nicefrac}
\usepackage{amsthm}
\usepackage{amsmath,amsfonts,bm,amssymb}
\usepackage{wrapfig}
\usepackage{caption}
\usepackage{siunitx}
\usepackage{thm-restate}
\usepackage{nccmath}
\usepackage{empheq}
\usepackage{bbm}
\usepackage{tabularx}

\usepackage{algorithmic}
\usepackage{algorithm}
\usepackage{filecontents}

\usepackage[capitalize,noabbrev]{cleveref}

\usepackage{color}
\usepackage{colortbl}
\definecolor{bgcolor}{rgb}{0.76,0.88,0.50}
\definecolor{bgcolor0}{rgb}{0.93,0.99,1}
\definecolor{bgcolor1}{rgb}{0.8,1,1}
\definecolor{bgcolor2}{rgb}{0.8,1,0.8}
\definecolor{bgcolor3}{rgb}{0.50,0.90,0.50}
\usepackage{tcolorbox}
\usepackage{pifont}

\definecolor{mydarkgreen}{RGB}{39,130,67}
\definecolor{mydarkorange}{RGB}{236,147,14}
\definecolor{mydarkred}{RGB}{192,47,25}
\definecolor{ruby}{RGB}{155,17,30}
\definecolor{chili}{RGB}{191,0,0}
\definecolor{sangria}{RGB}{146,0,10}
\definecolor{burgundy}{RGB}{128,0,32} 
\definecolor{darkred}{RGB}{132,0,0} 
\definecolor{cherry}{RGB}{192,0,0} 

\definecolor{blue}{RGB}{0,0,255}

\newcommand{\red}{\color{cherry}}

\usepackage[textsize=tiny]{todonotes}

\usepackage{xspace}
\newcommand{\algname}[1]{{\sf\footnotesize\red#1}\xspace}

\newcommand{\ind}{\perp\!\!\!\!\perp}
\newcommand{\norm}[1]{\left\| #1 \right\|}
\newcommand{\sqnorm}[1]{\left\| #1 \right\|^2}

\newcommand{\R}{\mathbb{R}} 
\newcommand{\N}{\mathbb{N}} 
\newcommand{\E}[1]{\mathbb{E}\left[#1\right]}

\newcommand{\Var}[1]{\mathrm{Var}\left(#1\right)}

\newcommand{\Exp}[1]{{\mathbb{E}}\left[#1\right]}
\newcommand{\ExpSub}[2]{{\mathbb{E}}_{#1}\left[#2\right]}

\newcommand{\cT}{\mathcal{T}}

\newcommand{\cJ}{\mathcal{J}}

\theoremstyle{plain}
\newtheorem{theorem}{Theorem}[section]
\newtheorem{proposition}[theorem]{Proposition}
\newtheorem{lemma}[theorem]{Lemma}

\theoremstyle{definition}

\newtheorem{assumption}[theorem]{Assumption}
\theoremstyle{remark}
\newtheorem{remark}[theorem]{Remark}

\newcommand{\eqdef}{:=}

\makeatletter
\newcommand{\vast}{\bBigg@{4}}

\DeclareMathOperator*{\argmin}{arg\,min}

\def\<{\left\langle}
\def\>{\right\rangle}
\def\[{\left[}
\def\]{\right]}
\def\({\left(}
\def\){\right)}

\usepackage{thmtools}
\usepackage{thm-restate}

\theoremstyle{theorem}
\newtheorem{innercustomthm}{Theorem}
\newenvironment{restate-theorem}[1]
{\renewcommand\theinnercustomthm{#1}\innercustomthm}
{\endinnercustomthm}

\newtheorem{innercustomlemma}{Lemma}
\newenvironment{restate-lemma}[1]
{\renewcommand\theinnercustomlemma{#1}\innercustomlemma}
{\endinnercustomlemma}

\newtheorem{innercustomproposition}{Proposition}
\newenvironment{restate-proposition}[1]
{\renewcommand\theinnercustomproposition{#1}\innercustomproposition}
{\endinnercustomproposition}

\newcommand*{\sketchproofname}{Sketch of Proof}

\usepackage{longtable}

\title{
  MindFlayer SGD: Efficient Parallel SGD in the Presence of\\
  Heterogeneous and Random Worker Compute Times
}

\author[1]{Artavazd Maranjyan}
\author[1]{Omar Shaikh Omar}
\author[1]{Peter Richt\'{a}rik}

\affil[1]{%
    King Abdullah University of Science and Technology\\
    Thuwal, Saudi Arabia
}
  
\begin{document}

\maketitle

\begin{abstract}

We investigate the problem of minimizing the expectation of smooth nonconvex functions in a distributed setting with multiple parallel workers that are able to compute stochastic gradients.
A significant challenge in this context is the presence of arbitrarily heterogeneous and stochastic compute times among workers, which can severely degrade the performance of existing parallel stochastic gradient descent (\algname{SGD}) methods.
While some parallel \algname{SGD} algorithms achieve optimal performance under deterministic but heterogeneous delays, their effectiveness diminishes when compute times are random—a scenario not explicitly addressed in their design.
To bridge this gap, we introduce \algname{MindFlayer SGD}, a novel parallel \algname{SGD} method specifically designed to handle stochastic and heterogeneous compute times.
Through theoretical analysis and empirical evaluation, we demonstrate that \algname{MindFlayer SGD} consistently outperforms existing baselines, particularly in environments with heavy-tailed noise.
Our results highlight its robustness and scalability, making it a compelling choice for large-scale distributed learning tasks.
\end{abstract}

\section{Introduction}
\label{section_intro}

We address the nonconvex optimization problem:
\begin{align}
  \label{eq_main_task}
  \min_{x \in \R^d} \Big \{f(x) \eqdef \ExpSub{\xi \sim \mathcal{D}}{f(x;\xi)}\Big \},
\end{align}
where $f\,:\,\R^d \times \mathbb{S} \rightarrow \R,$ and $\xi$ is a random variable with some distribution $\mathcal{D}$ on
$\mathbb{S}$.
In the context of machine learning, $\mathbb{S}$ could represent the space of all possible data, $\mathcal{D}$ denotes the distribution of the training dataset, and $f(\cdot, \xi)$ denotes the loss of a data sample $\xi$. 

The function $f$ is assumed to be differentiable, and its gradient is $L$--Lipschitz continuous (see Assumptions~\ref{ass:lipschitz_constant}--\ref{ass:lower_bound}).
We assume that we have $n$ workers available to work in parallel, each able to compute independent, unbiased stochastic gradients of $f$, whose variance is bounded by $\sigma^2$ (see Assumption \ref{ass:stochastic_variance_bounded}).
In this paper, we study the time complexity of methods working in this setup. 


We also assume access to $n$ parallel workers capable of computing independent stochastic gradients, in which case the classical approach is \algname{Minibatch SGD} \citep{cotter2011better,goyal2017accurate,gower2019sgd}.

\subsection{Minibatch SGD}
\label{section:minibatch_sgd}

\algname{Minibatch SGD} awaits the completion of all workers' computations of a single stochastic gradient before executing a gradient-type step:
\begin{enumerate}
    \item receive a single stochastic gradient $\nabla f(x^k; \xi^k_i)$ from each worker $i \in [n]$,
    \item update the model:
        $$
        x^{k+1} = x^{k} - \frac{\gamma}{n} \sum_{i=1}^n \nabla f(x^k; \xi^k_i),
        $$
\end{enumerate}
where $[n] \eqdef \{1,\ldots,n\}$, $\gamma>0$ is a stepsize, $\xi^k_i$ are i.i.d.\ samples from $\mathcal{D},$ and the gradients $\nabla f(x^k; \xi^k_i)$ are calculated in parallel.

In real systems, each worker's computational power may differ from the others, leading to varying completion times of gradient computation.
A notable drawback of \algname{Minibatch SGD} is its failure to account for these differences in compute times across workers.
The duration of each step is determined by the slowest worker's computation time. 
As a result, all other workers remain idle after completing their tasks, waiting for the slowest device to finish.
Meanwhile, this idle time could potentially be used in a more efficient way to improve the overall time complexity.
Clearly, a redesign of the algorithm is necessary.

\subsection{Asynchronous SGD} 
As a result, a new generation of algorithms emerged—asynchronous stochastic gradient descent (\algname{ASGD}) methods—designed to fully utilize all available computational resources \citep{recht2011hogwild,agarwal2011distributed,feyzmahdavian2016asynchronous,mania2017perturbed,nguyen2018sgd,arjevani2020tight,cohen2021asynchronous,mishchenko2022asynchronous,koloskova2022sharper,islamov2023asgrad,maranjyan2025ringmaster}.

Here, the server performs a gradient-type update immediately after receiving a stochastic gradient from any worker, without waiting for the others. 
The updated model is then sent back to the worker, which immediately begins computing a new stochastic gradient based on the updated model.
By the time the worker finishes computing this gradient, the model may have already been updated multiple times on the server due to gradients received from other workers. 
This creates a delay in the model update, denoted as $\delta_k$.
The algorithm can be described as follows:
\begin{enumerate}
    \item receive a stochastic gradient $\nabla f(x^{k-\delta_k}; \xi^{k-\delta_k})$ from any worker,
    \item update the model: 
        $$
        x^{k+1} = x^{k} - \gamma \nabla f(x^{k-\delta_k}; \xi^{k-\delta_k}),
        $$
    \item send new $x^{k+1}$ to the worker so the worker computes $\nabla f(x^{k+1}; \xi^{k+1})$.
\end{enumerate}
\citet{cohen2021asynchronous,mishchenko2022asynchronous,koloskova2022sharper} showed that \algname{ASGD} is provably faster in terms of time complexity than \algname{Minibatch SGD}.

However, it turns out that this untamed and wild asynchrony can be detrimental.
The drawback of \algname{ASGD} lies in the assumption that all worker computations are beneficial.
It suffers from the issue of updating the model with potentially significantly delayed gradients, which ultimately harms convergence and, consequently, the overall time complexity, as discussed in the work of \citet{tyurin2024optimal}.
To address this, there was a need to introduce a method that ignores outdated gradients while still maximizing the utilization of available computational resources.

\subsection{Rennala SGD}
Such a method was proposed in a recent breakthrough by \citet{tyurin2024optimal}.
Their method, \algname{Rennala SGD}, is a semi-asynchronous variant of \algname{Minibatch SGD}.
At each iteration, the server asynchronously collects a batch of gradients, allowing workers to send as many gradients as they can on the same point $x^k$.
Then, using this batch, \algname{Rennala SGD} proceeds with a gradient-type update using this batch as in \algname{Minibatch SGD}:
\begin{enumerate}
    \item wait until the server receives $B$ stochastic gradients at point $x^k$,
    \item update the model: 
        $$
        x^{k+1} = x^{k} - \frac{\gamma}{B} \sum_{j=1}^B \nabla f(x^k; \xi^k_j),
        $$
\end{enumerate}
more details on \algname{Rennala SGD} are in \Cref{sec_rennala}.
In this case, the faster the worker, the more gradients it sends.
For the struggling workers, it may happen that they are ignored.

\paragraph*{Worker Time.}
Their approach assumes each worker $i$ requires a fixed $\tau_i > 0$ seconds to compute a stochastic gradient.
For the first time lower bounds on time complexity were obtained for first order asynchronous methods in the fixed compute time regime for nonconvex functions with Lipschitz gradients.
They showed that \algname{Rennala SGD} is mini-max optimal in this setup in terms of time complexity.

Although it may seem like the story ends here, we challenge the fixed time assumption, arguing that a random time model better reflects reality.
The claim of optimality no longer holds due to this randomness, suggesting that the algorithms should be reevaluated and redesigned.
In this paper, we focus on this redesign, aiming to better align the algorithms with a more realistic model.

\section{Problem Setup and Contributions}


The deterministic compute time setup considered by \citet{tyurin2024optimal}, where \algname{Rennala SGD} is optimal, fails to capture the complexities of real-world distributed learning environments.
In practice, compute times are often uncertain due to various factors such as failing hardware, preemption by other jobs, delays in GPU computation, and inconsistencies in network communications \citep{chen2016revisiting, dutta2018slow}.
This uncertainty is even more pronounced in federated learning scenarios, where client unreliability can lead to unpredictable computation times or even incomplete tasks \citep{kairouz2021advances}.

To address these real-world challenges, we propose a practical setup that incorporates randomness in compute times.
Specifically, we consider a scenario where the stochastic gradient computation time of worker $i$ is given by:
\begin{equation}
    \label{eq:times}
    \tau_i +\eta_i,
\end{equation}
where $\tau_i > 0$ is a constant representing the minimum time for client $i$ to complete the gradient computation, and $\eta_i$ is a non-negative random variable drawn from some distribution $\mathcal{J}_i$, modeling the aforementioned uncertainties.

In this more realistic setting, existing methods like \algname{Rennala SGD} and \algname{ASGD} can perform poorly or even fail to converge.
We can illustrate this with a simple example:

\begin{figure*}[t]
	\centering
        \includegraphics[width=0.98\textwidth]{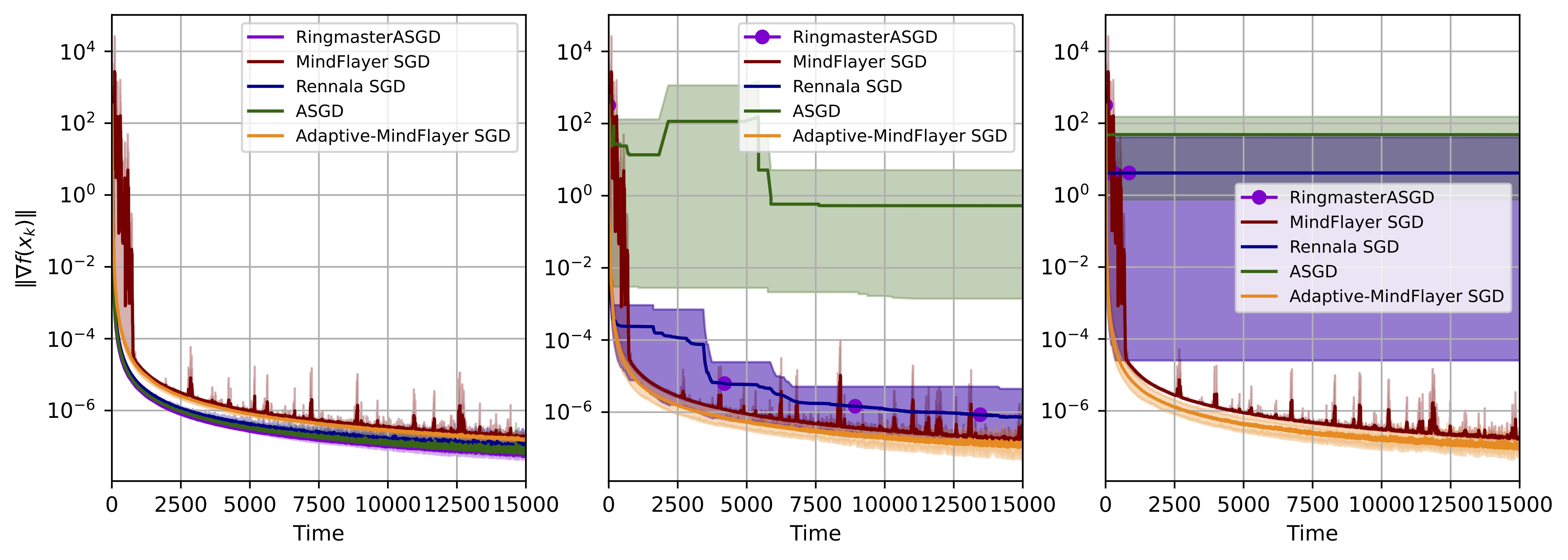}
	\caption{
    We conducted an empirical evaluation on a quadratic optimization problem (see \Cref{section_experiments} for details on the setup).
    In this experiment, we modeled the computation time for each worker as $\tau_i + \eta_i$, where $\tau_i = \sqrt{i}$ is a fixed constant and $\eta_i \sim \text{lognormal}(0, s)$ introduces noise for each worker $i \in [n]$.
    To investigate the effect of increasing noise variance, we tested different values of $s$: $s=1$ (left plot), $s=10$ (middle plot), and $s=100$ (right plot).
    As the variance increases, \algname{MindFlayer SGD} remains robust to the noise in worker computation times, whereas the convergence of the other methods \algname{Rennala SGD}, \algname{ASGD}, and \algname{Ringmaster ASGD}\protect\footnotemark{}, degrade significantly.
    This highlights how \algname{MindFlayer SGD} maintains efficiency even in the presence of heavy-tailed noise.
    Finally, we observe that \algname{Adaptive-MindFlayer SGD}, which does not assume prior knowledge of the distributions and adjusts the thresholds for each client adaptively based on observed times, performs better than \algname{MindFlayer SGD}. This improvement comes from its ability to adaptively select batch sizes for each worker based on random times, using a single batch size that is the sum of all individual ones.
    }
	\label{fig_modified_mf}
\end{figure*}
\footnotetext{This algorithm is presented by \citet{maranjyan2025ringmaster} and was published after our work. It is included here due to its relevance, but despite its novelty, our method still outperforms \algname{Ringmaster ASGD}.}

Consider a scenario where each time we request a device to compute a stochastic gradient, one of two outcomes occurs.
Either the device completes the computation exactly after the minimum time $\tau$ without any delays, or something goes wrong and the computation is never completed.
This situation can be modeled using a random time $\eta$ as follows:
\begin{equation}\label{eq:0_infinity}
  \eta = 
    \begin{cases} 
        0, & \text{with probability } 1-q, \\
        \infty\footnotemark, & \text{with probability } q,
    \end{cases}
\footnotetext{We can view $\eta$ as an extended real random variable, or just assume that $\infty$ is a very big number.}
\end{equation} 
where $0<q<1$.

In this scenario, any method that waits for a certain number of batches on each iteration to perform a step runs the risk of never receiving the required batch and getting stuck. 
This includes methods like \algname{Rennala SGD} or \algname{ASGD}. 
Specifically, if the algorithm waits for a single stochastic gradient on each iteration, with probability $q^n$, it will never receive it and consequently never proceed.

To address these limitations, we propose a new method that, unlike \algname{Rennala SGD} or \algname{ASGD}, does not wait for a fixed number of gradients.
Instead, it allocates a specific time for computing each stochastic gradient.
If a client fails to complete its computation within the designated time, the partial computation is discarded, and a new computation is initiated.
Our main contributions are as follows.
\begin{itemize}
    \item 
    In \Cref{section:mindflayer}, we introduce a time-efficient parallel \algname{SGD} method, \algname{MindFlayer SGD} (\Cref{alg_mind_flayer}), designed for the heterogeneous and random worker compute time regime described in \eqref{eq:times}.
    To the best of our knowledge, \algname{MindFlayer SGD} is the first algorithm tailored to this regime. We demonstrate that our method generalizes \algname{Rennala SGD}, making it optimal in the deterministic compute time setup.
    Furthermore, when the distribution of computation times is positively skewed, we show that \algname{MindFlayer SGD} outperforms the other methods, with the performance gap widening as the skewness coefficient increases.
    As illustrated in \Cref{fig_modified_mf}, where $\mathcal{J}_i = \text{lognormal}(0, s)$, increasing $s$ results in a higher skewness coefficient, which worsens the performance of \algname{Rennala SGD} and \algname{ASGD}.
    In contrast, \algname{MindFlayer SGD} remains robust to changes in variance.

    Additional experiments using various functions and distributions are presented in \Cref{section_experiments}.
    For distributions, we explore lognormal, log-Cauchy, and Infinite-Bernoulli (defined in \eqref{eq:0_infinity}).
    For functions, we examine a quadratic loss and a neural network trained on the \href{https://yann.lecun.com/exdb/mnist/}{MNIST} \citep{lecun1998gradient} dataset.
    This diverse testing framework highlights \algname{MindFlayer SGD}'s robustness and effectiveness across a wide range of challenging scenarios.
        
    \item 
    In \Cref{section:practicalmindflayer}, we introduce \algname{Adaptive-MindFlayer SGD}, a version of our algorithm that enhances practicality in two ways: it adapts to the computation time distribution during the learning process, eliminating the need for prior knowledge of the distribution, and it reduces the number of hyperparameters by treating all workers as a single entity, using only two hyperparameters in total.
    These improvements make \algname{Adaptive-MindFlayer SGD} more suitable for real-world implementation.

    \item 
    In \Cref{sec:comparison}, we compare our theoretical time complexity with that of \algname{Rennala SGD}, which is optimal in the deterministic time setting.
    We show that the time complexity of \algname{Rennala SGD} worsens as the distribution's tails become heavier or as the skewness coefficient increases, leading to an arbitrary performance gap compared to \algname{MindFlayer SGD}.

    \item 
    In \Cref{section:heteromindflayer}, we expand our theory to develop \algname{Vecna SGD}, designed for the heterogeneous case, where workers have datasets that are coming from different distributions. 
\end{itemize}

\section{Related Work}
\label{section_related_work}

There are several other related works. 
\citet{dutta2018slow} explore the error-runtime trade-offs in distributed SGD, revealing how slower and stale gradients can sometimes enhance convergence processes. 
\citet{woodworth2020local} compare local SGD with minibatch SGD, analyzing the efficiency of local updates in different distributed settings. 
\citet{wu2022delay} advance the understanding of asynchronous methods by proposing delay-adaptive step-sizes that adjust to asynchronous learning environments, optimizing the convergence rates. 
Furthermore, \citet{hanna2022adaptive, hanna2020strag} focus on adaptive stochastic gradient descent to improve communication efficiency in distributed learning, offering strategies that reduce communication demands while maintaining fast convergence. 

\section{Motivation and Single Device Case}
\label{section:single_device}

To illustrate the motivation behind the design of our method, let us consider a single device setup.
Recall the scenario introduced in \cref{eq:0_infinity} where we have single device and it either returns a gradient after $\tau$ time or gets stuck with probability $q$.
A straightforward and optimal workaround to this issue is to wait exactly $\tau$ seconds.
If we do not receive a gradient within this time frame, it indicates that we will never receive it, so there is no point in waiting longer.
In this case, we discard the current computation, which would never finish anyway, and request the device to compute the gradient again.
The probability of getting stuck again is lower, so eventually, we will receive a gradient and proceed.

Consider the following two strategies for each step.
\begin{itemize}
	\item \textbf{Strategy 1:} \algname{Rennala SGD}.
    We wait for the first $B$ stochastic gradients. 
    Thus, the time for one step for this strategy is the random variable:
    $$
      T_B = \sum_{j = 1}^B (\tau + \eta^j).
    $$
	\item \textbf{Strategy 2:} \algname{MindFlayer SGD}.
    We repeat the following random trial $B$ times: allocate time $\tau + t$ to computing a stochastic gradient.
    If no gradient is received within that time, discard the current computation and start over.
    Then the time for the $j$-th trial is given by:
		\begin{equation*}
			T^j(t) = 
			\begin{cases}
				\tau + \eta^j, & \text{if}\enspace \eta^j \le t, \\
				\tau + t,      & \text{if}\enspace \eta^j > t.
			\end{cases}
		\end{equation*}
    Thus, the time per step is a random variable:
    $$
        \tilde{T}_B(t) = \sum_{j=1}^{B} T^j(t).
    $$
\end{itemize}

\begin{figure*}[t]
	\centering
	\includegraphics[width=0.98\textwidth]{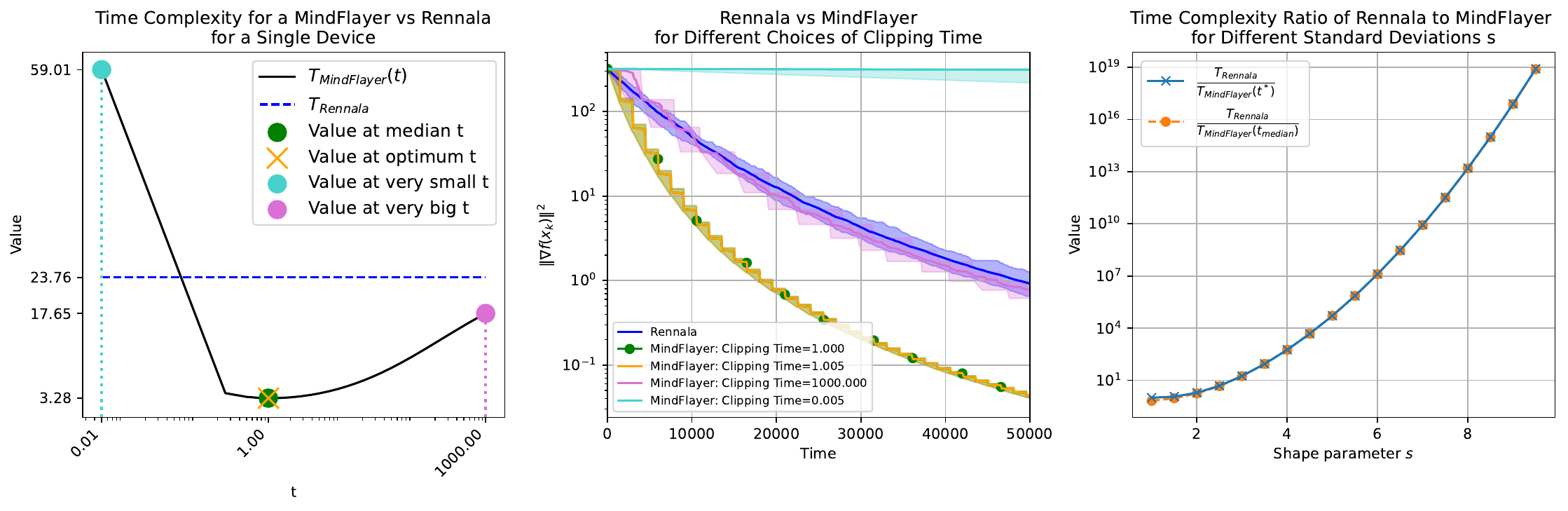}
	\caption{
        We consider a single-device setting with computation time $\tau + \eta$, where $\tau = 1$ and $\eta \sim \mathrm{lognormal}(0, s)$ for $s = 2.5$.
        \textbf{On the left}, we compare the expected time complexity of \algname{MindFlayer SGD} \eqref{eq:time_mindflayer} as a function of the time threshold $t$ with the constant time complexity of \algname{Rennala SGD} \eqref{eq:time_rennala}, showing the efficiency of \algname{MindFlayer SGD} for various choices of $t$.
        We set $K = 1$ since it appears in both algorithms and does not affect the comparison.
        Four key values of $t$ are highlighted: a very small one, the optimal choice from \eqref{eq:best_t}, the median of $\eta$, and a very large number.
        As expected, \algname{MindFlayer SGD} degrades for small $t$ due to the server receiving fewer gradients, remains unchanged between the median and optimal values, and worsens as $t$ increases, approaching \algname{Rennala SGD}.
        \textbf{In the middle}, we present an empirical evaluation on a quadratic optimization problem (see \Cref{section_experiments}), confirming time complexity reductions for \algname{MindFlayer SGD} across the same four thresholds.
        \textbf{On the right}, we plot the ratio of time complexities between \algname{Rennala SGD} and \algname{MindFlayer SGD} across different standard deviations $s$, revealing exponential efficiency gains at optimal clipping times, with similar trends at median values.
    }
	\label{fig_single_device_plots}
\end{figure*}

In the second case, rather than waiting for $B$ gradients, we attempt to compute $B$ gradients. 
Essentially, we limit the time spent on computing a stochastic gradient. 
In expectation, Strategy 2 will collect $Bp$ gradients per iteration, where $p = P(\eta \leq t)$ is the probability of collecting a gradient within a trial. 
Setting $t=\infty$ removes this restriction, resulting in the same strategy as the first one.

The fact that \algname{MindFlayer SGD} receives $Bp$ gradients on average makes it effectively a scaled-down version of \algname{Rennala SGD}.
Consequently, \algname{MindFlayer SGD} is expected to require $1/p$ times more iterations than \algname{Rennala SGD} to achieve the same level of convergence.
However, this trade-off is intentional and necessary to handle the heavy-tailed nature of computation times.

We have the following proposition.
\begin{proposition}[Proof in \Cref{section:proof_propositions}]
    \label{proposition:iterations}
    Consider the single device setup.
    Let $K$ be the number of iterations required by \algname{Rennala SGD} with batch size $B$ to find an $\varepsilon$-stationary point.
    For sufficiently small $\varepsilon$, \algname{MindFlayer SGD} with the same batch size $B$, needs $K/p$ iterations to find an $\varepsilon$-stationary point.
\end{proposition}

Thus, the time complexities in this setting are given by:
\begin{eqnarray}
	&& \label{eq:time_rennala}
  T_\algname{Rennala SGD} 
  = K \E{T_B} = K B(\tau + \E{\eta}),\\
	&& \notag
  T_\algname{MindFlayer SGD}(t)
  = \frac{K}{p} \E{\tilde{T}_B(t)} \\
  &&  \notag
  \quad = \frac{K}{p}B(\tau + (1-p)t + p\E{\tau | \tau \le t}) \\
  && \ \label{eq:time_mindflayer}
  \quad \leq \frac{K}{p}B(\tau + t).
\end{eqnarray}
This leads us to the following remark.
\begin{remark}
    \label{ass:t}
    For the case where $n=1$, \algname{MindFlayer SGD} is faster than \algname{Rennala SGD} if there exists a time threshold $t > 0$ such that the following inequality holds:
    $$
      P(\eta \le t) \ge \frac{\tau+t}{\tau + \E{\eta}}.
    $$
    The optimal choice of $t$ is given by
    \begin{equation} \label{eq:best_t}
        t = \argmin_{t > 0} \frac{\tau + t}{P(\eta \le t)}.
    \end{equation}
\end{remark}
Note that it is necessary that $t\le\E{\eta}$.
In case where $\E{\eta} = \infty$, this assumption hold for any finite values of $t>0$.
An example of this scenario is given in \eqref{eq:0_infinity}.
Many other distributions also have an infinite expectation, including Pareto, log-Cauchy, Lévy, log-t, and Landau distributions, among others.

A less restrictive example of distributions are positively skewed distributions.
Let 
$$
s = \E{\eta} - \mathrm{Med}[\eta]
$$ 
be the skewness coefficient of the distribution $\cJ$, where $\mathrm{Med}[\eta]$ denotes the median of $\eta$.
If $s>0$ we say that the distribution is positively skewed.
Then we have the following proposition.

\begin{proposition}[Proof in \Cref{section:proof_propositions}]
    \label{proposition:median}
    For the $n=1$ case, if $s>\tau + \mathrm{Med}[\eta]$ then \algname{MindFlayer SGD} is faster than \algname{Rennala SGD}.
    Moreover, if $s=\(\tau + \mathrm{Med}[\eta] \)\(2\alpha - 1\)$ then
    $$
      \frac{T_\algname{Rennala SGD}}{T_\algname{MindFlayer SGD}\(\mathrm{Med}[\eta]\)} \ge \alpha.
    $$
\end{proposition}
Therefore, \algname{Rennala SGD} can be arbitrarily bad.
As an example consider the lognormal$(\mu,\sigma^2)$ distribution. 
For this distribution, we have:
$$
s = \E{\eta} - \mathrm{Med}[\eta] = \exp\left(\mu + \frac{\sigma^2}{2}\right) - \exp(\mu).
$$

Thus, as we increase $\sigma$, the difference becomes arbitrarily large. 
To verify this, we also conducted a small experiment, see \Cref{fig_single_device_plots}. 
The right plot showcases how the ratio of time complexity between \algname{Rennala SGD} and \algname{MindFlayer SGD} can get arbitrarily large for the optimal clipping time 
$$
t^* := \arg \min_{t} T_{\algname{MindFlayer SGD}}(t)
$$
and even the median of the distribution $t_{\text{median}} = \mathrm{Med}[\eta]$. 
The left and middle plots showcase the potential improvement, and even loss resulting from different choices of clipping time $t$.

\section{MindFlayer SGD}
\label{section:mindflayer}

Here, we introduce \algname{MindFlayer SGD} for the multi-device setting ($n > 1$) \Cref{alg_mind_flayer}.
For the heterogeneous case, see \Cref{section:heteromindflayer}.

\begin{algorithm}[ht]
  \caption{\algname{MindFlayer SGD}\protect\footnotemark{}}
  \label{alg_mind_flayer}
  \begin{algorithmic}[1]
  \STATE \textbf{Input:} Initial point $x^0 \in \R^d$, stepsize $\gamma > 0$, allotted times $t_1, \dots, t_n \geq 0$, number of trials per client $B_1, \dots, B_n \geq 0$, probabilities $p_i = P(\eta_i \leq t_i)$
  \FOR{$k = 1, 2, \dots, K$}
      \STATE Send $x^k$ to all clients $i \in [n]$, each executes $B_i$ trials
      \STATE Compute the gradient estimate: \\
      \phantom{x}
      $
      g^k = \sum_{i=1}^{n} \sum_{j=1}^{B_i} I(\eta_i^j < t_i) \nabla f(x^k; \xi_i^j)
      $\label{asdf}
      \STATE Update: $x^{k+1} = x^{k} - \frac{\gamma}{B} g^{k}$, where $B=\sum_{i=1}^n p_i B_i$
  \ENDFOR
  \end{algorithmic}
\end{algorithm}

\algname{MindFlayer SGD} starts at an initial point $x^0 \in \R^d$ with stepsize $\gamma > 0$, time allowances $t_i > 0$, and trial counts $B_i \geq 0$ for each client.
At each iteration $k$, ranging from $1$ to $K$, the server sends the current point $x^k$ to all clients.

Each client $i$ makes $B_i$ attempts to compute stochastic gradients.
During each attempt, it computes a stochastic gradient, but if the computation exceeds the allotted time $\tau_i + t_i$, the gradient is discarded, and a new attempt begins.
This behavior is captured by the indicator function $I(\eta_i^j < t_i)$ in line \ref{asdf} of \Cref{alg_mind_flayer}.

\footnotetext{We name our method \algname{MindFlayer SGD}, drawing inspiration from \href{https://strangerthings.fandom.com/wiki/The_Mind_Flayer}{The Mind Flayer} from \textit{Stranger Things}, due to its ability to precisely control its clients, analogous to the creature's supreme control over its victims (The Flayed).}

The probability of completing the computation within the time limit is defined as $p_i \eqdef P(\eta_i^j < t_i)$.
Thus, the number of stochastic gradients received from client $i$ is a random variable ranging from $0$ to $B_i$, with an expected value of $p_i B_i$.
Summing over all clients, the expected total number of stochastic gradients is $B = \sum_{i=1}^n p_i B_i$.
Finally, after aggregating the collected gradients, the server updates the point using an unbiased gradient estimator, following the update rule $x^{k+1} = x^k - \frac{\gamma}{B} g^k$.

In the special case where the computation time is deterministic, i.e., $\eta_i = 0$ for every worker $i \in [n]$, we have $p_i = 1$ for all $i$. 
While \algname{Rennala SGD} does not explicitly specify the number of gradient computations $B_i$ for each client, in the deterministic setting, each client will send a fixed number of gradients per communication round. 
Consequently, for any $t > 0$, \algname{MindFlayer SGD} \Cref{alg_mind_flayer}, by choosing $B_i$ appropriately, reduces to \algname{Rennala SGD} \Cref{alg:rennala}.

However, the situation changes when $\eta_i > 0$ is not deterministic. 
If we set $t_i = \infty$ for all $i \in [n]$, \algname{MindFlayer SGD} \Cref{alg_mind_flayer} does not reduce to \algname{Rennala SGD} \Cref{alg:rennala}. 
This is because, in the case of \algname{Rennala SGD}, the randomness in each iteration causes the number of stochastic gradients computed by each client to vary across different communication rounds. 
Nevertheless, this scenario is not our primary focus, as we will demonstrate that allowing each worker to complete its gradient computation by setting $t_i = \infty$ is inefficient when dealing with positively skewed distributions.

Since \Cref{alg_mind_flayer} has multiple hyperparameters, we provide a more practical version in \Cref{section:practicalmindflayer}, where a single batch size $B$ and a single probability $p$ are used instead of client-specific values.
While this simplification is useful in practice, the more fine-grained version with per-client hyperparameters is generally more effective when computation time distributions are known.
In such cases, choosing the hyperparameters optimally is straightforward, making this version preferable when efficiency is a priority.
For this reason, we begin with a theoretical analysis of the more flexible version before introducing its practical counterpart.

\subsection{Theoretical Analysis}

We consider standard assumptions used in nonconvex optimization.


\begin{assumption}
    \label{ass:lipschitz_constant}
    Function $f$ is differentiable, and its gradient is $L$--Lipschitz continuous, i.e., 
    $$
        \norm{\nabla f(x) - \nabla f(y)} \leq L \norm{x - y}, \text{ for all } x, y \in \R^d.
    $$
 \end{assumption}

 \begin{assumption}
    \label{ass:lower_bound}
    The function $f(x)$ is bounded below, and we denote its infimum by $f^{\inf} \in \R$.
    Let $x^0$ be the initial point of the optimization method, define $\Delta \eqdef f(x^0) - f^{\inf}$.
\end{assumption}

\begin{assumption}
    \label{ass:stochastic_variance_bounded}
    For all $x \in \R^d,$ stochastic gradients $\nabla f(x;\xi)$ are unbiased and $\sigma^2$-variance-bounded, i.e., 
    \begin{align*}
        & \ExpSub{\xi}{\nabla f(x;\xi)} = \nabla f(x), \\
        & \ExpSub{\xi}{\norm{\nabla f(x;\xi) - \nabla f(x)}^2} \leq \sigma^2,
    \end{align*}
    where $\sigma^2 \geq 0$.
\end{assumption} 


\subsubsection{Convergence Theory}
The following theorem gives iterations guarantees for the convergence of \algname{MindFlayer SGD}.

Even though \algname{MindFlayer SGD} is similar to \algname{Rennala SGD} the convergence analysis require additional considerations, since the batch size is a random variable here as apposed to the case of \algname{Rennala SGD}. 

\begin{theorem}
  [Proof in \Cref{proof_theorem_mindflayer_convergence}]
    \label{theorem_mindflayer_convergence}
    Assume that Assumptions~\ref{ass:lipschitz_constant}, \ref{ass:lower_bound} and \ref{ass:stochastic_variance_bounded} hold. 
    Let     
    $$
        B = \sum_{i=1}^n p_i B_i \quad \text{and}\quad \gamma = \frac{1}{2L}\min\left\{1, \frac{\varepsilon B}{\sigma^2}\right\}
    $$
    in \Cref{alg_mind_flayer}.
    Then, the method guarantees that $\frac{1}{K}\sum_{k=0}^{K-1}\Exp{\norm{\nabla f(x^k)}^2} \leq \varepsilon$ after
    $$
        K \geq \max\left\{ 1, \frac{\sigma^2}{\varepsilon B} \right\} \frac{8 L \(f(x^0) - f^{\inf}\)}{\varepsilon}
    $$
    iterations.
\end{theorem}

Note that the rate is inversely proportional to the probabilities $p_i$, which is expected—smaller $p_i$ implies a potentially smaller batch size, leading to more iterations.

In the deterministic case where $\eta_i = 0$ for all $i \in [n]$, we have $p_i = P(\eta_i \le t_i) = 1$ for all $i \in [n]$, resulting in the same rate with $B = \sum_{i=1}^n B_i$.
This recovers the rate of \algname{Rennala SGD}, up to a constant factor.
Similarly, as $t_i \to \infty$ for all $i$, we get $p_i \to 1$, again leading to the same rate.

On the other hand, if $t_i = 0$ for all $i \in [n]$, then $K = \infty$, which is expected—if the success probability is zero for all clients, the server never receives any stochastic gradients, making progress impossible.


\subsubsection{Time Complexity}
The following theorem gives time complexity for \algname{MindFlayer SGD}.




\begin{theorem}[Proof in \Cref{proof_theorem_mindflayer_time_complexity}]
    \label{theorem_mindflayer_time_complexity}
    Assume that Assumptions~\ref{ass:lipschitz_constant}, \ref{ass:lower_bound} and \ref{ass:stochastic_variance_bounded} hold. 
    Let 
    $$
        B = \sum_{i=1}^n p_i B_i \quad \text{and}\quad \gamma = \frac{1}{2L}\min\left\{1, \frac{\varepsilon B}{\sigma^2}\right\}
    $$
    in \Cref{alg_mind_flayer}.
    Let $t = \(t_1,\ldots,t_n\)$, $t_1,\ldots,t_n \ge0$.
    Without loss of generality assume that 
    $$
        0< \tau_1 + t_1 \le \cdots \le \tau_n + t_n.
    $$
    Let 
    $$
        t(m) = \(\sum_{j=1}^{m}\frac{p_j}{\tau_j + t_j}\)^{-1}\(S+\sum_{j=1}^{m}p_j\),
    $$
    where $S = \max\left\{ 1, \frac{\sigma^2}{\varepsilon} \right\}$.
    Let $m^* = \arg\min_{m\in [n]} t(m)$, if there are several minimizers we take the smallest one.
    Put 
    $$
        B_i =  \lceil b_i \rceil, 
    \quad     
        b_i = 
        \begin{cases}
            \frac{t(m^*)}{\tau_i + t_i} - 1, & \text{if} \enspace i\le m^*,\\
            0, & \text{if} \enspace i > m^*.
        \end{cases}
    $$
    Then, \algname{MindFlayer SGD} guarantees to find an $\epsilon$-stationary point within
    \begin{align*}      
      &
        T_{\algname{MindFlayer SGD}}(t) \\
      &
        \le \frac{8 \Delta L}{\varepsilon} \times \min_{m\in [n]} \left\{ \(\frac{1}{m}\sum_{j=1}^{m}\frac{p_j}{\tau_j + t_j}\)^{-1}\(\frac{S}{m} + \bar{p}\) \right\}
    \end{align*}
    seconds, where $\bar{p} = \frac{1}{m}\sum_{j=1}^{m}p_j$.
\end{theorem}

The term $\nicefrac{\tau_i + t_i}{p_i}$ also appears in the single-device case (\Cref{ass:t}), where the optimal choice of $t_i$ minimizes this quantity, similar to \eqref{eq:best_t}.
As in the single-device case, we observe an inverse dependence on $p_i$: smaller $p_i$ increases $\nicefrac{\tau_i + t_i}{p_i}$, which may lead to certain devices being excluded from $m^*$.
This is expected, as a small $p_i$ indicates an unreliable device.

For the first $m^*$ selected workers, the optimal allocation of $B_i$ depends on $\tau_i + t_i$: the smaller this value (i.e., the faster the device), the more trials it should receive.
Thus, the choice of $B_i$ takes into account both the speed of a device ($\tau_i$) and its reliability ($p_i$).
Some devices may be fast (small $\tau_i$) but unreliable (small $p_i$), leading to $B_i = 0$ for such devices.

This allocation strategy ensures that devices with high computation times and low reliability are excluded.
For example, a device might be computationally fast but suffer from frequent network issues, preventing gradients from reaching the server.
In such cases, the server should not rely on this device and should exclude it from gradient computations.

In the deterministic case where $\eta_i = 0$, we have $p_i = 1$ for all $i \in [n]$.
The optimal choice in this setting is $t_i = 0$, which recovers the same time complexity as \algname{Rennala SGD}.

This theorem assumes that the computation time distributions of the devices are known.
While this may seem restrictive, it is often reasonable in distributed systems, where models are trained repeatedly.
Over time, these repeated executions provide enough data to approximate the distributions with high confidence, making such an assumption practical in many cases.

However, relying on this prior knowledge is not always feasible, especially in settings like federated learning, where devices may be unpredictable or newly introduced into the system.
To address this, in the next section, we present a practical version of the algorithm that does not require knowledge of computation time distributions, making it more adaptable to real-world scenarios.

\section{Practical MindFlayer SGD}
\label{section:practicalmindflayer}

To run our algorithm, we need to specify $B_i$ for each client.
One possible approach is to learn the distribution of each client's behavior on the fly and assign an appropriate $B_i$ accordingly.
Similar ideas are explored by \citet{maranjyan2025ata}.
However, this is not straightforward in our setting, since we perform thresholding and cannot directly observe or learn the distribution of clients' compute times.
Instead, in this work, we replace the client-specific parameters $B_i$ with a single global parameter $B$.
To make this possible, we treat all clients in a uniform way.
If every client had the same probability $p$ of successfully sending a gradient of equal quality, such a simplification would be justified.

To enforce this uniform probability, we need to set $t_i$ such that each client has a probability $p$ of completing its computation within the given threshold.
Since we do not assume prior knowledge of the computation time distribution, we estimate $t_i$ empirically.
For each client, we need to solve the following problem: given $p$, find $t$ such that
$$
    P(\tau + \eta \le t) = p.
$$
Since $P(\tau + \eta < t)$ is a non-decreasing function of $t$ and $I(\tau + \eta \leq t)$ is an unbiased estimator, we can solve this using Robbins-Monro stochastic approximation method \citep{robbins1951stochastic}.
This algorithm iteratively updates the threshold $t$ using the rule:
\begin{equation}\label{eq:robbins-monro}
    t^{r+1} = t^r - \alpha^r \( I(\tau + \eta^r \leq t^r) - p \),
\end{equation}
where $\alpha^r$ is the stepsize, which we take $\alpha^r = \alpha^0/i$ with $\alpha^0 > 0$.

Note that we do not need to know $\tau$ and $\eta^r$; we only require $I(\tau + \eta^r \leq t^r)$, which is 1 if the worker finishes the computation within the threshold and 0 otherwise.

Putting all the pieces together we derive our algorithm \Cref{alg:mod_mf}.

\begin{algorithm}[ht] 
    \caption{\algname{Adaptive-MindFlayer SGD}}
    \label{alg:mod_mf}
    \begin{algorithmic}[1]
        \STATE \textbf{Input:} 
        Initial point $x^0 \in \R^d$, stepsize $\gamma > 0$, probability $p \in (0, 1]$, number of trials $B \geq 0$
        \STATE \textbf{Robbins-Monro inputs:} initial thresholds $t_i^0$, initial stepsizes $\alpha_i^0>0$ for all $i\in [n]$
        \FOR{$k = 1, 2, \dots, K$}
          \STATE Each worker $i \in [n]$ computes a gradient at $x^k$
          \STATE Initialize $g^k = 0$ and $b = 0$
          \WHILE{$b < B$} 
            \STATE Wait for the fastest client $i$ to finish its trial
            \STATE Receive $g = I(\tau_i + \eta_i^{k_b} \le t_i^{k_b}) \nabla f(x^k; \xi_i^{k_b})$
            \STATE Update $t_i^{k_b}$ using \eqref{eq:robbins-monro}
            \STATE Update: $g^k = g^k + g$; $b = b + 1$
          \ENDWHILE
          \STATE Update $x^{k+1} = x^k - \frac{\gamma}{pB} g^{k}$ 
        \ENDFOR
    \end{algorithmic}
\end{algorithm}

The algorithm is controlled by two key parameters: $p$ and $B$.
The parameter $p$ represents system reliability and is shared across all clients, allowing us to define a single global parameter $B$.
Although the Robbins-Monro process requires per-client parameters, these are not critical for the algorithm's performance.
For the initial threshold, one can simply start with a large value, as the algorithm will automatically adjust it over time.
Similarly, the step size is not highly sensitive; in our experiments, we used $1$, but the algorithm performs well with other choices.

By dynamically adjusting $t_i$ based on real-time observations of worker compute times, \algname{Adaptive-MindFlayer SGD} continuously aligns the clipping threshold with the desired completion probability $p$.
This eliminates the need for extensive manual tuning of hyperparameters and improves robustness to variability in compute times.

Even without assuming prior knowledge of computation time distributions, this empirical threshold selection performs nearly as well as the previous approach, which relied on such knowledge.
In \Cref{fig_modified_mf}, we show that \algname{Adaptive-MindFlayer SGD} achieves comparable performance to \algname{MindFlayer SGD} while simplifying hyperparameter selection, making it particularly practical for distributed systems with heterogeneous and unpredictable worker compute times.

\section{Comparing to Rennala SGD}
\label{sec:comparison}

\begin{figure*}[t]
	\centering
	\includegraphics[width=0.98\linewidth]{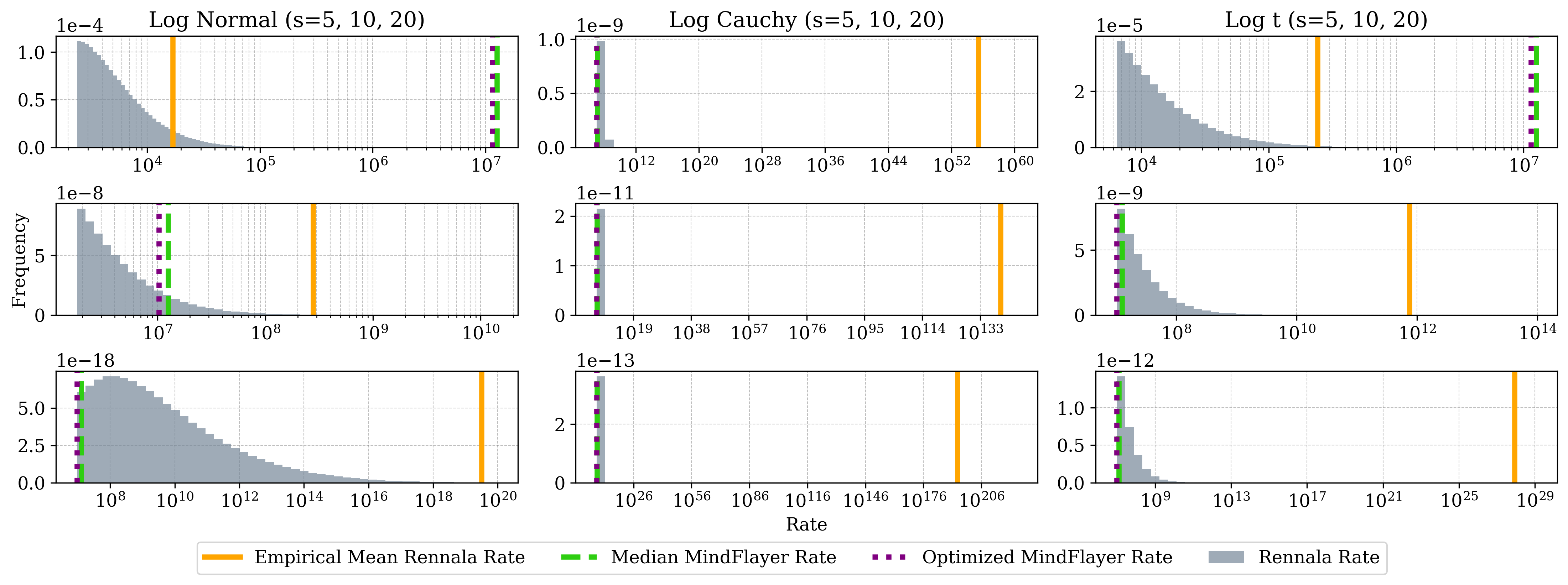}
	\caption{
        Empirical comparison of the theoretical time complexities of \algname{Rennala SGD} and \algname{MindFlayer SGD}.
        We set $n = 100$ and $B = 250$ (see \Cref{sec:exp_theoretical} for details).
        Since \algname{Rennala SGD}'s time complexity is a random variable, we plot a histogram with its empirical mean.
        For \algname{MindFlayer SGD}, we compare the theoretical time complexity from \Cref{theorem_mindflayer_time_complexity}, evaluating two strategies for selecting $t_i$:
        (1) using the median of the distributions $\mathcal{J}_i$, and (2) finding the optimal $t_i$ from \eqref{eq:best_t}.
        We examine three distributions: lognormal, log-Cauchy, and log-$t$ (5 degrees of freedom), shown in each column.
        As variance increases (rows: 5, 10, 20), \algname{MindFlayer SGD} increasingly outperforms \algname{Rennala SGD}.
    }
	\label{fig_distirbution_comparison}
\end{figure*}

Comparing the theoretical performance of \algname{Rennala SGD} and \algname{MindFlayer SGD} is particularly challenging due to the inherent randomness in the time complexity of \algname{Rennala SGD} and the dependence of \algname{MindFlayer SGD} on optimizing time variables $t_i$.
A comparison based on expected time complexity overlook the nuances of each algorithm's performance across different distributions.
Therefore, we turn to an empirical comparison to provide insights into their practical behavior.
Particularly, we demonstrate how \algname{MindFlayer SGD} can achieve arbitrarily small time complexity in heavy-tailed distributions.

To begin, we derive the time complexity of \algname{Rennala SGD} in the context of random times.
Let 
$$
    \mathcal{B}:=\left\{(B_1,B_2,\ldots,B_n): B_i\in \N_0; \sum_{i=1}^n B_i=B \right\}
$$
be the set of all possible batch sizes for each device, the time $T_B$ required for one step with batch size $B$ of \algname{Rennala SGD} is given by:
\begin{align*}
    T_B 
    &= \min_{\mathcal{B}} \left\{\max_{i\in[1,n]}\left\{B_i\tau_i + \sum_{j=1}^{B_i}\eta_i^j\right\}\right\} \ge T_1 \\
    & = \min_{i \in [n]} \left\{\tau_i + \eta_i^1\right\} \ge \min_{i \in [n]} \tau_i + \min_{i \in [n]}\eta_i^1.
\end{align*}
Thus, the expected time to collect a batch $B$ is
\begin{equation*}
    \E{T_B} \ge \tau_{\min} + \E{\min_{i \in [n]}\eta_i}.
\end{equation*}
Note that if the distribution of $\min_{i \in [n]}\eta_i$ is heavy-tailed, then the expected time complexity may be infinite, thus favoring \algname{MindFlayer SGD} over \algname{Rennala SGD}.
A simple illustration of this occurs when extending the \Cref{eq:0_infinity} case, where $\eta$ is either zero or infinite, to scenarios involving multiple devices.
In such cases, the expectation of the minimum time across devices, $\min_{i \in [n]}\eta_i$, also results in an infinite expected time complexity.

While a detailed theoretical comparison is intractable, we conduct an empirical comparison to highlight practical differences between the two algorithms.
To capture the randomness in \algname{Rennala SGD}'s rate, we generate a histogram of $T_B$ and convolve it with itself $K$ times, where $K$ is the number of iterations needed for $\epsilon$-convergence.
We set $n = 100$ and $B = 250$. Additional details are provided in \Cref{sec:exp_theoretical}.

For \algname{MindFlayer SGD}, we compare the theoretical time complexity from \Cref{theorem_mindflayer_time_complexity}.
We evaluate two strategies for selecting $t_i$:
(1) using the median of the distributions $\mathcal{J}_i$, and (2) solving the
following optimization problem: 
\begin{quote}
	Fix $m \in [n]$, minimize $t(m)$ over \\ 
    $t=(t_1,\cdots,t_n)$, (remember $p_j=F_j(t_j)$).
\end{quote}
We optimize this using the \texttt{L-BFGS-B} algorithm, a well-suited method
for solving smooth, convex, or mildly nonconvex problems due to its efficiency
and robustness \citep{zhu1997algorithm}.
For each $m$, we take the minimum over all possible configurations.

Our empirical results, illustrated in \Cref{fig_distirbution_comparison},
demonstrate that as the variance of the underlying distribution increases,
\algname{MindFlayer SGD} consistently outperforms \algname{Rennala SGD}.
The heavy-tailed nature of the distributions causes \algname{Rennala SGD} to
experience extreme slowdowns, while \algname{MindFlayer SGD} maintains robust
performance.

\section{Conclusion and Future Work} 
\label{section_conclusion_future_work}

In this paper, we address the problem of minimizing the expectation of nonconvex functions with Lipschitz gradients, with the use of parallel workers computing stochastic gradients. 
Our focus lies on the challenging scenario where worker compute times are heterogeneous and random, expanding on recent developments in \algname{ASGD} methods like \algname{Rennala SGD}. 
We observe that while \algname{Rennala SGD} performs optimally in environments with deterministic compute times, its effectiveness diminishes under random compute conditions.

To better understand and improve stochastic optimization in these conditions, we introduce a novel asynchronous \algname{SGD} method named \algname{MindFlayer SGD}. 
This method adjusts to the randomness in computation times by not adhering to a fixed batch size but rather setting specific times for computing single stochastic gradients. 
If a client fails to deliver within this time frame, the computation is discarded, and the process restarts. 
This flexibility allows \algname{MindFlayer SGD} to perform robustly across various conditions, notably outperforming both \algname{Rennala SGD} and standard Asynchronous \algname{SGD} (\algname{ASGD}) in our theoretical and empirical analysis.

Our results demonstrate that \algname{MindFlayer SGD} significantly reduces time complexity, particularly in environments characterized by positively skewed distribution of computation times. 
We empirically validate this in simulations with several distributions conditions where \algname{MindFlayer SGD} consistently outperforms the other methods, particularly in high-variance scenarios. 
This showcases its superiority in adapting to the unpredictable duration of gradient computations typical in real-world applications such as federated learning environments.

In this study, our analysis was confined to computation times, with no consideration given to communication times. 
Future research will extend our investigation to include communication times. 
Moreover, we plan to explore the application of gradient estimators with varying variance bounds across different clients. 
We hypothesize that controlling these variance bounds could yield further benefits in the optimization process.

\begin{acknowledgements} 
    The research reported in this publication was supported by funding from King Abdullah University of Science and Technology (KAUST): i) KAUST Baseline Research Scheme, ii) Center of Excellence for Generative AI, under award number 5940, iii) SDAIA-KAUST Center of Excellence in Artificial Intelligence and Data Science.
\end{acknowledgements}

\bibliography{bib}

\newpage
\onecolumn
\appendix




\section{Experimental details}
\label{section_experiments}

In this section we explain the setup for comparing \algname{MindFlayer SGD}, \algname{Rennala SGD}, and \algname{ASGD}, which we used throughout this paper. 
We compare the algorithms' performance on a quadratic optimization \eqref{eq:quad-problem} task with access to a stochastic gradient. 
The parallelism was simulated on a machine with 2 Intel(R) Xeon(R) Gold 6226R CPUs @ 2.90GHz, with a total of 64 logical CPUs. 
For each setting of the algorithm, we run 10 different seeds for the random time and plot the average, minimum and maximum, see \Cref{fig_modified_mf}, \Cref{fig_single_device_plots}, etc.

We use a similar setup to the one employed by \citet{tyurin2024optimal}, but modify it so that we have a known expected variance. 
We make this choice, so we can compare theoretical parameters, as we did in \Cref{fig_single_device_plots}.

Furthermore, we consider the homogeneous optimization problem \ref{eq_main_task}, with the
convex quadratic function:
\begin{equation}
\label{eq:quad-problem}
  f(x) = \frac{1}{2} x^\top A x - b^\top x \qquad \forall x \in \mathbb{R}^d~.
\end{equation}

We take $d = 1000$,


\begin{align*}
  \textstyle
  A = \frac{1}{4}
  \begin{bmatrix}
  2 & -1 &  & 0 \\
  -1 & \ddots & \ddots &  \\
  & \ddots & \ddots & -1 \\
  0 & & -1 & 2 \\
  \end{bmatrix}
  \in \mathbb{R}^{d \times d}
  \quad \text{and} \quad
  b = \frac{1}{4}
  \begin{bmatrix}
  -1 \\
  0 \\
  \vdots \\
  0 \\
  \end{bmatrix}
  \in \mathbb{R}^d.
\end{align*}

Assume that all $n$ workers has access to the following unbiased stochastic
gradients:
$$
\textstyle
\left[\nabla f(x, \xi)\right]_j := \nabla_j f(x) + \xi~,
$$
where $\xi \sim \mathcal{N}(0, 0.0003^2)$, thus, we get that in
\Cref{ass:stochastic_variance_bounded} we have, 
$$
\textstyle
\sigma^2 = 0.0003^2 \cdot d = 0.0003^2 \cdot 1000~.
$$

\begin{figure*}[t]
	\centering
 	\includegraphics[width=0.98\linewidth]{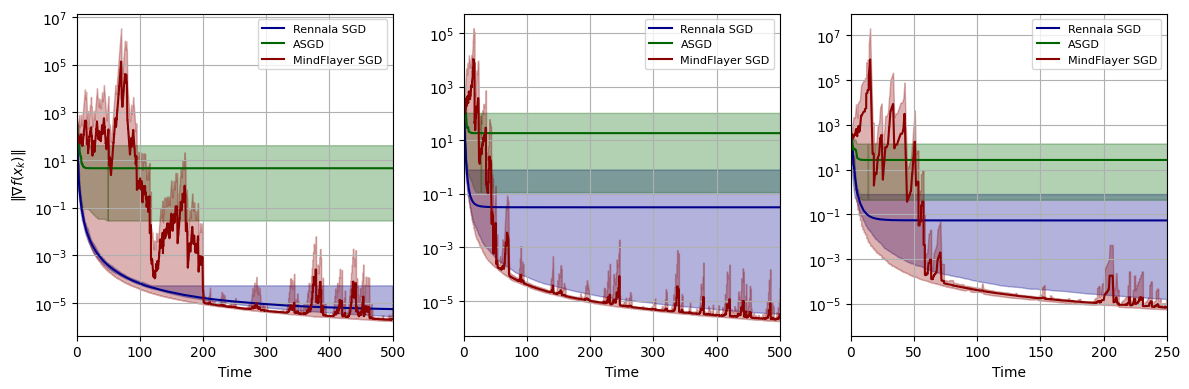}
	\caption{
        We ran an experiment as described in \Cref{section_experiments} where we employ the same $\mathcal{J}_i = \text{InfBernoulli}(q)$ distribution for all clients $i \in [n]$, with different $q$ values.
        From left to right we have $q = 0.6, 0.7, 0.8$.
        Additionally, we set $\tau_i = \sqrt{i}$.
        As we observe, with an increase of the probability of failure $q$ unlike \algname{Rennala SGD} and \algname{ASGD}, \algname{MindFlayer SGD} demonstrates the ability to continue optimizing and not be stuck
    }
 	\label{inf_bernoulli_multiple}
 \end{figure*}

\begin{figure*}[t]
	\centering
	\includegraphics[width=0.98\linewidth]{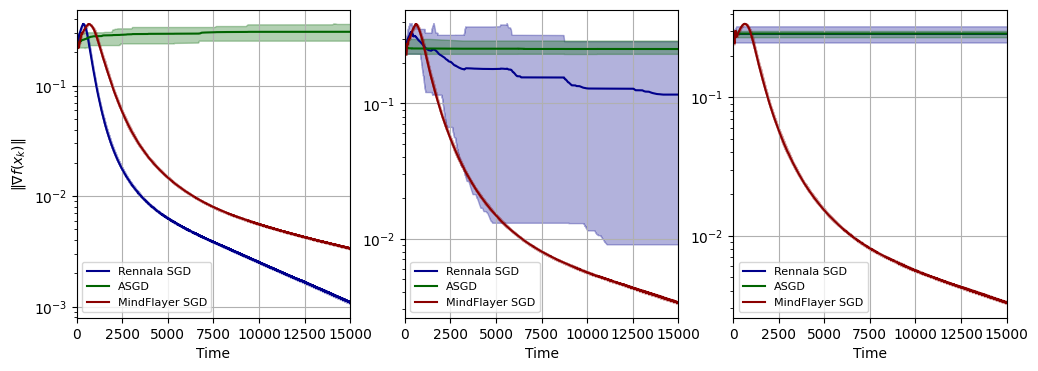}
	\caption{We train a two layer Neural Network on the MNIST dataset where we
	set the distribution $\mathcal{J}_i = \text{log-Cauchy}(s)$ for all clients
	$i \in [n]$, with different scale values $s$. From left to right we have $s
	= 1, 10, 100$.
    Additionally, we set $\tau_i = \sqrt{i}$.
    We observe that
	\algname{MindFlayer SGD} convergence doesn't suffer from the increase in
	the scale parameter $s$. On the other hand, \algname{Rennala} and \algname{ASGD}
	are delayed significantly with bigger scale parameters $s$}
	\label{log_cauchy}
\end{figure*}

Now setting the convergence threshold $\epsilon = 10^{-4}$, we can infer all theoretical parameters.
To find the optimal time corresponding to \algname{Rennala SGD} we need to fix
the times, we do that by either removing the randomness, or adding the expected
randomness. On the other hand, for \algname{MindFlayer SGD} we use the results
from \Cref{theorem_mindflayer_time_complexity} to set the theoretical number of
trials for each client. For some experiments we used theoretical stepsizes, 
e.g. Figure \ref{fig_single_device_plots}, for others we used the range of stepsizes
from a set $\{2^i | i \in [-10, 10]\}$, e.g. Figures
\ref{fig_modified_mf}, \ref{fig_modified_mf}, and
\ref{inf_bernoulli_multiple}, similarly to \citet{tyurin2024optimal}. Finally,
for the nonconvex problem in Figure \ref{log_cauchy} we tried the set $\{0.01,
0.001, 0.0001\}$.

In addition to the experimental results shown throughout the paper, we ran two
more experiments. One with the Infinite-Bernoulli distribution on the same 
quadratic problem, and a second with the log-Cauchy distribution with a
small two-layer neural network on the MNIST dataset, see \Cref{inf_bernoulli_multiple} and \Cref{log_cauchy}.

\subsection{Comparing Theoretical Time Complexities}
\label{sec:exp_theoretical}

To produce the empirical comparison demonstrated in
Figure~\ref{fig_distirbution_comparison}, we design a synthetic setup to
compare the time complexities of \algname{MindFlayer SGD} and \algname{Rennala
SGD}. We use distributions with varying tail behaviors and levels of variance.
Specifically, the base task time for each worker is defined as $ \tau_i = 100
\sqrt{i}$, which accentuates systematic variation in task durations across
workers. On top of this base time, random delays are added, drawn from one of
the distributions listed in Table~\ref{tab:random_times}.
These distributions
are parameterized by varying scale factors $ s = 5, 10, 20 $.
The log-$ t $
distribution uses $ \text{df} = 5 $ to control its tail behavior.

The rate of \algname{MindFlayer SGD} is derived from
Theorem~\ref{theorem_mindflayer_time_complexity}, expressed as:
$$
T_{\algname{MindFlayer SGD}}(t) \le 8 \times \min_{m \in [n]} \left\{ \left( \frac{1}{m} \sum_{j=1}^{m} \frac{p_j}{\tau_j + t_j} \right)^{-1} \left( \frac{S}{m} + \frac{1}{m} \sum_{j=1}^{m} p_j \right) \frac{\Delta L}{\varepsilon} \right\}
$$

Here, the rates depend on the clipping times $ t = (t_1, \dots, t_n) $, which
determine $ p_i \coloneqq F_i(t_i) $. As described in
Section~\ref{sec:comparison}, we use both the medians of the distributions and
the output of an optimizer as baselines for $ t $.

For \algname{Rennala SGD}, the time complexity is modeled as an approximation
of the random variable $ \sum_{i=1}^{K} T^i_B $, where $ K = \max\left\{1,
\frac{\sigma^2}{\varepsilon} \right\} \frac{8 \Delta L}{\varepsilon} $ is the
iteration complexity, as demonstrated in the proof of
Proposition~\ref{proposition:iterations} in
Section~\ref{section:proof_propositions}. The random variables $ T^i_B $
represent the time to collect a theoretical batch of size $ B = \lceil \sigma^2
/ \varepsilon \rceil $. As described in Section~\ref{sec:comparison}, we obtain
an approximation of the distribution of the sum using convolution.

Both rates depend on the problem constants $ \varepsilon, \sigma^2, L, \Delta
$, which are defined in Table~\ref{tab:parameters}.

\begin{table}[h!]
    \centering
    \begin{tabular}{|c|l|}
        \hline
        \textbf{Distribution} & \textbf{Parameters} \\ \hline
        log-Normal            & Scale $ s = 5, 10, 20 $, Mean $ \mu = 0 $ \\ \hline
        log-Cauchy            & Location $ \mu = 0 $, Scale $ s = 5, 10, 20 $ \\ \hline
        log-$ t $             & Degrees of freedom $ \text{df} = 5 $, Location $ \mu = 0 $, Scale $ s = 5, 10, 20 $ \\ \hline
    \end{tabular}
    \caption{Random Time Setup for Worker Task Times}
    \label{tab:random_times}
\end{table}

\begin{table}[h!]
    \centering
    \begin{tabular}{|c|l|c|}
        \hline
        \textbf{Symbol} & \textbf{Description} & \textbf{Value} \\ \hline
        $ \sigma $ & Standard deviation of noise & $ \sqrt{0.025} $ \\ \hline
        $ \epsilon $ & Convergence threshold & $ 10^{-4} $ \\ \hline
        $ \Delta \cdot L $ & Lipschitz constant multiplied by suboptimality & 1 \\ \hline
        $ B $ & \algname{Rennala SGD} batch size & $ \lceil \sigma^2 / \epsilon \rceil $ \\ \hline
        $ S $ & Theoretical Parameter in \algname{MindFlayer SGD} & $ \lceil \sigma^2 / \epsilon \rceil $ \\ \hline
        $ n $ & Number of workers & 100 \\ \hline
        $ n_{\text{samples}} $ & Number of random samples to estimate $T_B$ & 1000 \\ \hline
        $ \tau_i $ & Fixed compute time for worker $ i $ & $ 100 \sqrt{i} $ \\ \hline
    \end{tabular}
    \caption{Experimental Parameters}
    \label{tab:parameters}
\end{table}

\section{Table of Notations}

\begin{table}[H]
    \label{table:notation}
    \begin{center}
    \begin{tabular}{ll}
    \multicolumn{1}{c}{\bf Notation} & \multicolumn{1}{c}{\bf Meaning} \\ \hline \\
    $[n]$       & $\{1,\ldots, n\}$ \\
    $L$         & Lipschitz constant of gradients, i.e., $\norm{\nabla f(x) - \nabla f(y)} \leq L \norm{x - y}$ (\Cref{ass:lipschitz_constant}) \\
    $f^{\inf}$  & Minimum value of the function, i.e., $f^{\inf} \leq f(x)$ (\Cref{ass:lower_bound}) \\ 
    $\Delta$    & $\Delta \eqdef f(x^0) - f^{\inf}$, where $x^0$ is the starting point of optimization methods \\
    $\sigma^2$  & Variance bound on gradients, i.e., $\ExpSub{\xi}{\norm{\nabla f(x;\xi) - \nabla f(x)}^2} \leq \sigma^2$ (\Cref{ass:stochastic_variance_bounded}) \\
    $\gamma$    & Stepsize \\
    $\tau_i$    & Minimum time required for client $i$ to compute a gradient \\
    $\eta_i$    & Additional random time taken while computing the gradient \\
    $\cJ_i$     & Distribution of the non-negative random variable $\eta_i$ \\
    $t_i$       & Allotted time for worker $i$ to compute a gradient \\
    \end{tabular}
    \end{center}
\end{table}

\section{Heterogeneous Regime}
\label{section:heteromindflayer}

So far, we have discussed the scenario where all workers compute i.i.d. stochastic gradients. 
However, in distributed optimization and federated learning \citep{konevcny2016federated}, workers may have different datasets. 
Consider the following optimization problem:
\begin{align}
    \label{eq:main_task_heterog}
    \min_{x \in \R^d} \Big \{f(x) \eqdef \frac{1}{n} \sum_{i=1}^n \ExpSub{\xi_i \sim \mathcal{D}_i}{f_i(x;\xi_i)}\Big \},
\end{align}
where $f_i\,:\,\R^d \times \mathbb{S}_{i} \rightarrow \R^d$ and $\xi_i$ are random variables with some distributions $\mathcal{D}_i$ on $\mathbb{S}_{i}$.
Problem \eqref{eq:main_task_heterog} generalizes problem \eqref{eq_main_task}. 

\subsection{Related work and discussion}

The optimization problem \eqref{eq:main_task_heterog} has been thoroughly studied in many papers, including \citep{aytekin2016analysis, mishchenko2018delay, nguyen2022federated, wu2022delay, koloskova2022sharper, mishchenko2022asynchronous}.  
There have been attempts to analyze Asynchronous SGD in the heterogeneous setting.  
For example, \citet{mishchenko2022asynchronous} demonstrated convergence only to a neighborhood of the solution.  
In general, achieving good rates for Asynchronous SGD is difficult without making additional assumptions about the similarity of the functions $f_i$ \citep{koloskova2022sharper, mishchenko2022asynchronous}.

In the deterministic case, when $\sigma^2 = 0$, \citet{wu2022delay} analyzed the \algname{PIAG} method in the deterministic heterogeneous regime and showed convergence.  
Although the performance of \algname{PIAG} can be good in practice, in the worst case \algname{PIAG} requires $\operatorname{O}\left(\nicefrac{\tau_{n} \widehat{L} \Delta}{\varepsilon}\right)$ seconds to converge, where $\tau_{n}$ is the time delay of the slowest worker, $\widehat{L} \eqdef \sqrt{\sum_{i=1}^n L_i^2},$ and $L_i$ is a Lipschitz constant of $\nabla f_i$.  
Note that the synchronous \algname{Minibatch SGD} (see \Cref{section:minibatch_sgd}) method has the complexity $\operatorname{O}\left(\nicefrac{\tau_{n} L \Delta}{\varepsilon}\right)$, which is always better.\footnote{In the nonconvex case, $\widehat{L}$ can be arbitrarily larger than $L$.}  

\citet{tyurin2024optimal} proposed an optimal method in the regime where worker computation times are deterministic, similar to the homogeneous setup.

\subsection{Vecna SGD}

Here we describe our method called \algname{Vecna SGD}.

\begin{algorithm}[H]
    \caption{\algname{Vecna SGD}\protect\footnotemark{}}
    \label{alg_heteromindflayer}

    \begin{algorithmic}[1]
    \STATE \textbf{Input:} starting point $x^0 \in \R^d$, stepsize $\gamma>0$, allotted times $t_1,\ldots, t_n \ge 0$, number of trials per client $B_1,\ldots, B_n \ge 0$ 
    \FOR{$k = 1, 2, \dots, K$}
    \STATE Put $g_i^k = 0$
    \STATE Send $x^k$ to all clients
    \STATE Run Method~\ref{alg_vecna_client} in all clients $i=1,2,\ldots,n$
    \WHILE{there is a client that has trials to perform} 
    \STATE Wait for the fastest client
    \STATE Receive gradient $g_i$ from client $i$
    \STATE $g_i^k = g_i^k + g$
    \ENDWHILE
    \STATE $g^k = \frac{1}{n}\sum_{i=1}^n \frac{g_i^k}{p_iB_i}$, \hfil $\diamond$ $p_i = F_{i}(t_i) = P(\eta_i \le t_i)$.
    \STATE $x^{k+1} = x^k - \gamma g^{k}$ 
    \ENDFOR
    \end{algorithmic}
\end{algorithm}

\begin{algorithm}[H]
    \caption{Client $i$-s $k$-th step}
    \label{alg_vecna_client}

    \begin{algorithmic}[1]
    \STATE Receive $x^k$ from the server
    \FOR{$j = 1,2,\ldots, B_i$}
    \STATE Sample $\eta_i^j\sim \cJ_i$ \hfill $\diamond$ Start computing gradient estimator.
    \IF{$\eta_i^j \le t_i$} 
    \STATE $g = \nabla f(x^k;\xi_i^j), \enspace \xi_i^j \sim \mathcal{D}$ \hfill $\diamond$ The computation completes within the allotted time $t_i$.
    \STATE Send $g$ to the server
    \ENDIF
    \ENDFOR
    \end{algorithmic}
\end{algorithm}
\footnotetext{We name our method \algname{Vecna SGD}, drawing inspiration from \href{https://strangerthings.fandom.com/wiki/Vecna}{Vecna} from \textit{Stranger Things}.}

The \algname{Vecna SGD} algorithm begins with an initialization at a starting point $x^0$ in $\R^d$, with a specified stepsize $\gamma$, time allowances $t_i$, and trial counts $B_i$ for each client.
In each iteration $k$, ranging from $k = 1$ to $K$, the server distributes the current point $x^k$ to all clients.
Each client $i$ then executes a subroutine (\Cref{alg_vecna_client}) to attempt to compute $B_i$ stochastic gradients from samples $\xi_i^j$ drawn from a distribution $\mathcal{D}$.
During each attempt, client $i$ starts computing a stochastic gradient; if the computation exceeds the allotted time $t_i$, they discard the current gradient and begin another computation.
Consequently, the actual number of stochastic gradients received from each client $i$ becomes a random variable, ranging from 0 to $B_i$.
The expected number of gradients from client $i$ is given by $p_i B_i$.
The server normalizes the gradients by the expected batch size $p_i B_i$ and then aggregates them.
Finally, the point is updated to $x^{k+1} = x^k - \gamma g^k$ following each aggregation round.

\subsection{Convergence theory}

The following theorem gives iterations guarantees for the convergence of \algname{Vecna SGD}.

\begin{theorem}[Proof in \Cref{proof_heteromindflayer_convergence}]
  \label{theorem_heteromindflayer_convergence}
  Assume that Assumptions~\ref{ass:lipschitz_constant}, \ref{ass:lower_bound} and \ref{ass:stochastic_variance_bounded} hold for function $f_i$ for all $i\in [n]$.
  Let 
  $
  \gamma = \min\left\{\frac{1}{\sqrt{L\alpha K}}, \frac{1}{L \beta}, \frac{\varepsilon}{2 L \zeta}\right\}
  $
  in \Cref{alg_heteromindflayer}.
  Then after 
  \begin{align*}
      K \geq \frac{12 \Delta L}{\varepsilon}\max\left\{\beta, \frac{12\Delta \alpha}{\varepsilon}, \frac{2 \zeta}{\varepsilon}\right\},
  \end{align*}
  iterations, the method guarantees that $\min_{0\le k \le K} \E{\norm{\nabla f(x^k)}^2} \leq \varepsilon$, and 
  $$
    \alpha = \frac{L}{n} \max_{i\in [n]} \left\{\frac{1-p_i}{p_i B_i}\right\}, \quad 
    \beta = 1, \quad 
    \zeta = \frac{\sigma^2}{n^2}\sum_{i=1}^{n} \frac{1}{p_i B_i} 
        + \frac{2L}{n} \max_{i\in [n]} \left\{ \frac{1-p_i}{p_i B_i} \right\} \Delta^{\inf},
  $$
  where $\Delta^{\inf} \eqdef \sum_{i=1}^{n} \(f^{\inf} - f_i^{\inf}\)$.
\end{theorem}

\subsection{Time Complexity}

The following theorem gives time complexity for \algname{Vecna SGD}.

\begin{theorem}[Proof in \Cref{proof_theorem_heteromindflayer_time_complexity}]
  \label{theorem_heteromindflayer_time_complexity}
    Assume that Assumptions~\ref{ass:lipschitz_constant}, \ref{ass:lower_bound} and \ref{ass:stochastic_variance_bounded} hold for function $f_i$ for all $i\in [n]$.
  Let $\gamma = \min\left\{\frac{1}{\sqrt{L\alpha K}}, \frac{1}{L}, \frac{\varepsilon}{2 L }\right\}$ in \Cref{alg_heteromindflayer}, where
  $$
    \alpha = \frac{L}{n} \max_{i\in [n]} \left\{\frac{1-p_i}{p_i B_i}\right\}, \quad 
    \zeta = \frac{\sigma^2}{n^2}\sum_{i=1}^{n} \frac{1}{p_i B_i} 
        + \frac{2L}{n} \max_{i\in [n]} \left\{ \frac{1-p_i}{p_i B_i} \right\} \Delta^{\inf}.
  $$
  Let $t = \(t_1,\ldots,t_n\)$, $t_1,\ldots,t_n \ge0$.
  Without loss of generality assume that $0< \tau_1 + t_1 \le \cdots \le \tau_n + t_n$.
  Let 
  $$
    T = \tau_n + t_n + \[\frac{1}{n}\sum_{i=1}^n \frac{\tau_i + t_i}{p_i}\] \frac{\sigma^2}{n \varepsilon} 
        + \max_{i\in [n]} \left\{\frac{1-p_i}{p_i}\(\tau_i + t_i\)\right\} \frac{L\(\Delta + \Delta^{\inf}\)}{n\varepsilon}~.
  $$
  Put 
  $$
  B_i =  \lceil b_i \rceil, 
  \quad     
  b_i = \frac{T}{\tau_i + t_i}.
  $$
  Then, \algname{Vecna SGD} guarantees to find an $\epsilon$-stationary point after
  $$
  T_{\algname{Vecna SGD}}(t) 
  \le 288 \times \frac{\Delta L}{\varepsilon} 
  \(\tau_n + t_n + \[\frac{1}{n}\sum_{i=1}^n \frac{\tau_i + t_i}{p_i}\] \frac{\sigma^2}{n \varepsilon} 
        + \max_{i\in [n]} \left\{\frac{1-p_i}{p_i}\(\tau_i + t_i\)\right\} \frac{L\(\Delta + \Delta^{\inf}\)}{n\varepsilon}\)
  $$
  seconds. 
\end{theorem}

\section{The Rennala Algorithm}
\label{sec_rennala}

\begin{minipage}{0.55\textwidth}
\begin{algorithm}[H]
    \caption{\algname{Rennala SGD}}
    \label{alg:rennala}
    \begin{algorithmic}[1]
    \STATE \textbf{Input:} starting point $x^0$, stepsize $\gamma$, batch size $S$
    \STATE Run Method~\ref{alg:alg_worker} in all workers
    \FOR{$k = 0, 1, \dots, K - 1$}
    \STATE Init $g^k = 0$ and $s = 1$
    \WHILE{$s \leq S$}
    \STATE Wait for the next worker
    \STATE Receive gradient and iteration index $(g, k')$
    \IF{$k' = k$}
    \STATE $g^k = g^k + \frac{1}{S} g$; \;\; $s = s + 1$
    \ENDIF
    \STATE Send $(x^k, k)$ to the worker
    \ENDWHILE
    \STATE $x^{k+1} = x^k - \gamma g^{k}$
    \ENDFOR
    \end{algorithmic}
\end{algorithm}
\end{minipage}
\hfill
\begin{minipage}{0.45\textwidth}
\begin{algorithm}[H]
    \caption{Worker's Infinite Loop}
    \label{alg:alg_worker}
    \begin{algorithmic}[1]
    \STATE Init $g = 0$ and $k' = -1$
    \WHILE{True}
    \STATE Send $(g, k')$ to the server
    \STATE Receive $(x^k, k)$ from the server
    \STATE $k' = k$
    \STATE $g = \nabla f(x^k;\xi), \quad \xi \sim \mathcal{D}$
    \ENDWHILE
    \end{algorithmic}
\end{algorithm}
\end{minipage}

We mention the \algname{Rennala SGD} throughout the paper, here we provide a brief introduction to the method and its development. 
\Cref{alg:rennala} shows the work done by the server. 
Essentially, the server asynchronously waits to collect a batch of size $S$, whenever it receives a gradient from a worker that has the same iteration as the algorithm, it assigns it to compute a gradient at the same point $x_k$. 
After collecting the batch, we preform a synchronous update (given that all gradients were made on the same point $x_k$), using an average of the collected batch.

\section{The classical SGD theory}
In this section, we present the classical \algname{SGD} theory as developed by \citet{ghadimi2013stochastic} and \citet{khaled2020better}. 
Our analysis will follow the approach of the latter. 

We consider the stochastic gradient descent (\algname{SGD}) method:
$$
x^{k+1} = x^k - \gamma g(x^k),
$$
where $x^0\in \R^d$ is the initial point, and $g(x)$ is a stochastic gradient estimator at $x$.

We make the following assumption:
\begin{assumption}
    \label{ass:ABC}
    The stochastic gradient estimator $g(x)$ satisfies:
    \begin{align*}
        &\E{g(x)} = \nabla f(x)\\
        &\Exp{\norm{g(x)}^2} \leq 2\alpha\(f(x) - f^{\inf}\) + \beta \norm{\nabla f(x)}^2 + \zeta,
    \end{align*} 
    for all $x\in \R^d$ and some constants $\alpha, \beta, \zeta \ge 0$.
\end{assumption}

This assumption is both general and reasonable, and it is satisfied by many modern \algname{SGD}-type methods. 
For further details, refer to \citet{khaled2020better}.

Under this assumption, we can derive the following convergence result.

\begin{theorem}[Corollary 1 \citep{khaled2020better}]
    \label{theorem_sgd}
    Assume that Assumptions~\ref{ass:lipschitz_constant}, \ref{ass:lower_bound} and \ref{ass:ABC} hold. 
    Then for any $\varepsilon>0$
    \begin{align*}
        \min_{0\le k \le K} \E{\norm{\nabla f(x^k)}^2} \leq \varepsilon
    \end{align*}
    for
    \begin{align*}
        \gamma = \min\left\{\frac{1}{\sqrt{L\alpha K}}, \frac{1}{L \beta}, \frac{\varepsilon}{2 L \zeta}\right\},
    \end{align*} 
    and
    \begin{align*}
        K \geq \frac{12 L \(f(x_0) - f^{\inf}\)}{\varepsilon}\max\left\{\beta, \frac{12\Delta \alpha}{\varepsilon}, \frac{2 \zeta}{\varepsilon}\right\}.
    \end{align*}
\end{theorem}

\section{Proofs for Propositions in Section \ref{section:single_device}}
\label{section:proof_propositions}

\begin{restate-proposition}{\ref{proposition:iterations}}
  Consider the single device setup.
  Let $K$ be the number of iterations required by \algname{Rennala SGD} with batch size $B$ to find an $\varepsilon$-stationary point.
  For sufficiently small $\varepsilon$, \algname{MindFlayer SGD} with the same batch size $B$, needs $K/p$ iterations to find an $\varepsilon$-stationary point.
\end{restate-proposition}
\begin{proof}
    The iterations of \algname{Rennala SGD} can be viewed as iterations of \algname{Minibatch SGD}. 
    Thus, we can apply the classical \algname{SGD} theory (\Cref{theorem_sgd}) to derive its iteration complexity:
    $$
    K = \max\left\{ 1, \frac{\sigma^2}{\varepsilon B} \right\} \frac{8 L (f(x^0) - f^{\inf})}{\varepsilon}.
    $$    
    For \algname{MindFlayer SGD}, the iteration complexity follows from \Cref{theorem_mindflayer_convergence}. 
    Therefore, the number of iterations $K_M$ required for \algname{MindFlayer SGD} to guarantee that 
    $$
    \frac{1}{K_M}\sum_{k=0}^{K_M-1}\Exp{\norm{\nabla f(x^k)}^2} \leq \varepsilon
    $$
    is given by
    $$
    K_M = \max\left\{ 1, \frac{\sigma^2}{\varepsilon Bp} \right\} \frac{8 L (f(x^0) - f^{\inf})}{\varepsilon}.
    $$
    If $\varepsilon \leq \frac{\sigma^2}{B}$, we have
    $$
    K_M = \frac{K}{p}.
    $$
\end{proof}

\begin{restate-proposition}{\ref{proposition:median}}
    For the $n=1$ case, if $s>\tau + \mathrm{Med}[\eta]$ then \algname{MindFlayer SGD} is faster than \algname{Rennala SGD}.
    Moreover, if $s=\(\tau + \mathrm{Med}[\eta] \)\(2\alpha - 1\)$ then
    $$
    \frac{T_\algname{Rennala SGD}}{T_\algname{MindFlayer SGD}\(\mathrm{Med}[\eta]\)} \ge \alpha.
    $$
\end{restate-proposition}

\begin{proof}
    Let $t= \mathrm{Med}[\eta] =: m$, recall that $s = \E{\eta} - m$, then we have:
    $$
        T_\algname{MindFlayer SGD}(m) \leq \frac{K}{p}B(\tau + t) = 2KB\(\tau + m\),
    $$
    $$
        T_\algname{Rennala SGD} = K B(\tau + \E{\eta}) = K B(\tau + m + s).
    $$
    Thus, if $s>\tau + m$ then \algname{MindFlayer SGD} is faster than \algname{Rennala SGD}.

    Now, let $s=\(\tau + m \)\(2\alpha - 1\)$ then
    $$
    \frac{T_\algname{Rennala SGD}}{T_\algname{MindFlayer SGD}\(m\)} \ge \frac{\tau + m + s}{2\(\tau + m\)} = \frac{2\alpha\(\tau + m\)}{2\(\tau + m\)} = \alpha.
    $$
\end{proof}
\section{Proofs for Homogeneous Regime}

\subsection{Proof of Theorem~\ref{theorem_mindflayer_convergence}}
\label{proof_theorem_mindflayer_convergence}

First, we rewrite \algname{MindFlayer SGD} in a classical \algname{SGD} way where we do gradient step with an unbiased estimator of the gradient at each iteration.

\begin{algorithm}[H]
    \caption{\algname{MindFlayer SGD}}
    \label{alg_mind_flayer_sgd}
    \begin{algorithmic}[1]
    \STATE \textbf{Input:} starting point $x^0$, stepsize $\gamma$, time budgets $t_1,\ldots,t_n\ge 0$, batch sizes $B_1,\ldots,B_n \ge 0$,
    \FOR{$k = 0, 1, \dots, K - 1$}
    \STATE $g^k = \frac{1}{B}\sum_{i=1}^{n}\sum_{j=1}^{B_i}I\(\eta_i^j \le t_i\) \nabla f\(x^k;\xi_i^j\)$
    \STATE $x^{k+1} = x^k - \gamma g^{k}$
    \ENDFOR
    \end{algorithmic}
\end{algorithm}
where $B=\sum_{i=1}^n p_i B_i$, $p_i = F(t_i) = P(\eta_i \le t_i)$ and $I(\cdot)$ denotes the indicator function.
To prove the theorem we need to establish some properties of the gradient estimator.
First, we need an unbiased estimator.

\begin{lemma}[Proof in \Cref{proof_lemma_mindflayer_unbiased}]
    \label{lemma_mindflayer_unbiased}
    The gradient estimator in \Cref{alg_mind_flayer_sgd} given by 
    $$
    g(x) \eqdef \frac{1}{B}\sum_{i=1}^{n}\sum_{j=1}^{B_i}I\(\eta_i^j \le t_i\) \nabla f\(x;\xi_i^j\)
    $$
    is unbiased, i.e., $\E{g(x)} = \nabla f(x)$ for all $x\in\R^d$.
\end{lemma}

Next, we obtain an upper bound for the variance of this estimator.

\begin{lemma}[Proof in \Cref{proof_lemma_mindflayer_variance}]
    \label{lemma_mindflayer_variance}
    The gradient estimator in \Cref{alg_mind_flayer_sgd} given by 
    $$
    g(x) \eqdef \frac{1}{B}\sum_{i=1}^{n}\sum_{j=1}^{B_i}I\(\eta_i^j \le t_i\) \nabla f\(x;\xi_i^j\)
    $$
    satisfies
    $$
    \E{\norm{g(x)^2}} \leq 2\sqnorm{\nabla f(x)} + \frac{1}{B}\sigma^2.
    $$
\end{lemma}

We are ready to prove the \Cref{theorem_mindflayer_convergence}.

\begin{restate-theorem}{\ref{theorem_mindflayer_convergence}}
    Assume that Assumptions~\ref{ass:lipschitz_constant}, \ref{ass:lower_bound} and \ref{ass:stochastic_variance_bounded} hold. 
    Let 
    $$
        B = \sum_{i=1}^n p_i B_i \quad\text{and}\quad \gamma = \frac{1}{2L}\min\left\{1, \frac{\varepsilon B}{\sigma^2}\right\}
    $$
    in \Cref{alg_mind_flayer}.
    Then, after
    $$
    K \geq \max\left\{ 1, \frac{\sigma^2}{\varepsilon B} \right\} \frac{8 L \(f(x^0) - f^{\inf}\)}{\varepsilon}
    $$
    iterations, the method guarantees that $\frac{1}{K}\sum_{k=0}^{K-1}\Exp{\norm{\nabla f(x^k)}^2} \leq \varepsilon$.
\end{restate-theorem}

\begin{proof}
    Note that \Cref{alg_mind_flayer} can be viewed as a special case of classical stochastic gradient descent (SGD), as reformulated in \Cref{alg_mind_flayer_sgd}. 
    We need to verify that the gradient estimator fulfills the conditions required by classical \algname{SGD} (\Cref{theorem_sgd}). 
    The two preceding lemmas address this requirement precisely. 
    Specifically, \Cref{lemma_mindflayer_unbiased} confirms that the gradient estimator used in \Cref{alg_mind_flayer_sgd} is unbiased, while \Cref{lemma_mindflayer_variance} verifies that the variance of this estimator meets the conditions specified in \Cref{ass:ABC}, with $\alpha = 0$, $\beta = 2$ and $\zeta = \frac{\sigma^2}{B}$.
    Consequently, it remains to apply \Cref{theorem_sgd}.
\end{proof}

\subsubsection{Proof of \Cref{lemma_mindflayer_unbiased}}
\label{proof_lemma_mindflayer_unbiased}

\begin{restate-lemma}{\ref{lemma_mindflayer_unbiased}}
    The gradient estimator in \Cref{alg_mind_flayer_sgd} given by 
    $$
      g(x) := \frac{1}{B}\sum_{i=1}^{n}\sum_{j=1}^{B_i}I\(\eta_i^j \le t_i\) \nabla f\(x;\xi_i^j\)
    $$
    is unbiased, i.e., $\E{g(x)} = \nabla f(x)$ for all $x\in\R^d$,
    where $B=\sum_{i=1}^n p_i B_i$.
\end{restate-lemma}

\begin{proof}
    This follows from direct computation:
    \begin{eqnarray*}
        \E{g(x)} 
        &=& \E{\frac{1}{B}\sum_{i=1}^{n}\sum_{j=1}^{B_i}I\(\eta_i^j \le t_i\) \nabla f\(x;\xi_i^j\)}\\
        &=& \frac{1}{B}\sum_{i=1}^{n}\sum_{j=1}^{B_i}\E{I\(\eta_i^j \le t_i\) \nabla f\(x;\xi_i^j\)}\\
        &\overset{\(\eta_i^j \ind \xi_i^j\)}{=}& \frac{1}{B}\sum_{i=1}^{n}\sum_{j=1}^{B_i}\E{I\(\eta_i^j \le t_i\)} \E{\nabla f\(x;\xi_i^j\)}\\
        &=& \frac{1}{B}\sum_{i=1}^{n}\sum_{j=1}^{B_i}p_i \nabla f(x)\\
        &=& \nabla f(x) \frac{1}{B}\sum_{i=1}^{n} p_i B_i \\
        &=& \nabla f(x).
    \end{eqnarray*}
\end{proof}

\subsubsection{Proof of \Cref{lemma_mindflayer_variance}}
\label{proof_lemma_mindflayer_variance}

\begin{restate-lemma}{\ref{lemma_mindflayer_variance}}
    The gradient estimator in \Cref{alg_mind_flayer_sgd} given by 
    $$
    g(x) := \frac{1}{B}\sum_{i=1}^{n}\sum_{j=1}^{B_i}I\(\eta_i^j \le t_i\) \nabla f\(x;\xi_i^j\)
    $$
    satisfies
    $$
    \E{\norm{g(x)^2}} \leq 2\sqnorm{\nabla f(x)} + \frac{1}{B}\sigma^2,
    $$
    where $B=\sum_{i=1}^n p_i B_i$.
\end{restate-lemma}

\begin{proof}
    In order to simplify notation, let 
    $$
    a_i := \sum_{j=1}^{B_i}b_i^j,
    $$
    where
    $$
    b_i^j := I\(\eta_i^j \le t_i\) \nabla f\(x;\xi_i^j\).
    $$

\textbf{Step 1 (Initial expression).} 
We express $\E{\|g(x)\|^2}$ in terms of $a_i$:
\begin{eqnarray*} 
    \E{\sqnorm{g(x)}} 
    &=& \E{\sqnorm{\frac{1}{B}\sum_{i=1}^{n}a_i}} = \frac{1}{B^2}\E{\sum_{i=1}^{n}\sqnorm{a_i} + \sum_{i \ne j} \<a_i, a_j\>}.
\end{eqnarray*}

We further simplify both terms via:
\begin{eqnarray}\label{a_i}
    \sqnorm{a_i} = \sqnorm{\sum_{j=1}^{B_i} b_i^j} = \sum_{j=1}^{B_i} \sqnorm{b_i^j} + \sum_{k\ne l} \<b_i^k, b_i^l\>,
\end{eqnarray}

\begin{eqnarray}\label{a_ia_j}
    \<a_i, a_j\> = \<\sum_{k=1}^{B_i}b_i^k, \sum_{l=1}^{B_j}b_j^l \> = \sum_{k=1}^{B_i} \sum_{l=1}^{B_j} \<b_i^k, b_j^l\>.
\end{eqnarray}
\textbf{Step 2. (Finding the expectations).}
Further
    \begin{eqnarray}\label{b_i}
        \E{\sqnorm{b_i^j}} 
        &=& \E{\(I\(\eta_i^j \le t_i\)\)^2\sqnorm{\nabla f\(x;\xi_i^j\)}} \notag\\
        &\overset{\(\eta_i^j \ind \xi_i^j\)}{=}& \E{\(I\(\eta_i^j \le t_i\)\)^2}\E{\sqnorm{\nabla f\(x;\xi_i^j\)}} \notag\\
        &\le& p_i \(\sqnorm{\nabla f(x)} + \E{\sqnorm{\nabla f\(x;\xi_i^j\) - \nabla f(x)}}\) \notag\\
        &\overset{\(\Cref{ass:stochastic_variance_bounded}\)}{\le}& p_i \(\sqnorm{\nabla f(x)} + \sigma^2\),
    \end{eqnarray}
    and 
    \begin{eqnarray}\label{b_ib_j}
        \E{\<b_i^k, b_j^l\>} 
        &=& \E{\<I\(\eta_i^k \le t_i\) \nabla f\(x;\xi_i^k\), I\(\eta_j^l \le t_j\) \nabla f\(x;\xi_j^l\) \>} \notag \\ 
        &\overset{\( \ind\)}{=}& \E{I\(\eta_i^k \le t_i\)} \E{I\(\eta_j^l \le t_j\)} \<\E{ \nabla f\(x;\xi_i^k\)} , \E{\nabla f\(x;\xi_j^l\)} \> \notag \\
        &=& p_i p_j \sqnorm{\nabla f(x)}.
    \end{eqnarray}
\textbf{Step 3 (Putting everything together).} 
We start with
    \begin{eqnarray*}
    \E{\sqnorm{a_i}} 
    &\overset{\(\ref{a_i}, \ref{b_i}, \ref{b_ib_j}\)}{\le}& B_i p_i \(\sqnorm{\nabla f(x)} + \sigma^2\) + B_i\(B_i - 1\)p_i^2 \sqnorm{\nabla f(x)}\\
    &\le& B_i p_i \(\sqnorm{\nabla f(x)} + \sigma^2\) + B_i^2p_i^2 \sqnorm{\nabla f(x)},
    \end{eqnarray*}
using this and recalling the definition of $B$, we get
    \begin{eqnarray*}
        \E{\sum_{i=1}^{n}\sqnorm{a_i}}
        &\le& B \sqnorm{\nabla f(x)} + B \sigma^2 + \sqnorm{\nabla f(x)} \sum_{i=1}^{n}B_i^2p_i^2.
    \end{eqnarray*}
    Next
    \begin{eqnarray*}
        \<a_i, a_j\> \overset{\(\ref{a_ia_j}, \ref{b_ib_j}\)}{=} B_i p_i B_j p_j \sqnorm{\nabla f(x)},
    \end{eqnarray*}
    finally,
    \begin{eqnarray*}
        \E{\sqnorm{g(x)}} 
        &=& \frac{1}{B^2}\E{\sum_{i=1}^{n}\sqnorm{a_i} + \sum_{i \ne j} \<a_i, a_j\>}\\
        &\le& \frac{1}{B^2}\[ B \sqnorm{\nabla f(x)} + B \sigma^2 + \(\sum_{i=1}^{n}B_i^2p_i^2 + \sum_{i \ne j} B_i p_i B_j p_j \) \sqnorm{\nabla f(x)}\]\\
        &=& \frac{1}{B^2}\( B + B^2 \) \sqnorm{\nabla f(x)} + \frac{\sigma^2}{B}\\
        &\le&  2\sqnorm{\nabla f(x)} + \frac{\sigma^2}{B}.
    \end{eqnarray*}
\end{proof}

\subsection{Proof of Theorem~\ref{theorem_mindflayer_time_complexity}} 
\label{proof_theorem_mindflayer_time_complexity}

The following lemma gives time complexity for any choice of $B_1,\ldots, B_n$ and $t = \(t_1,\ldots,t_n\)$ in \algname{MindFlayer SGD}.
\begin{lemma}[Proof in \Cref{proof_lemma_mindflayer_time_complexity}]
    \label{lemma_mindflayer_time_complexity}
    Assume that Assumptions~\ref{ass:lipschitz_constant}, \ref{ass:lower_bound} and \ref{ass:stochastic_variance_bounded} hold.
    Let 
    $$
        B = \sum_{i=1}^n p_i B_i \quad\text{and}\quad \gamma = \frac{1}{2L}\min\left\{1, \frac{\varepsilon B}{\sigma^2}\right\}
    $$
    in \Cref{alg_mind_flayer}. 
    Then after 
    $$
        T_{\algname{MindFlayer SGD}}(t) 
        \le \max_{i \in [n]}\{ B_i\(\tau_i + t_i\) \} \max\left\{ 1, \frac{\sigma^2}{\varepsilon B} \right\} \frac{8 L \(f(x_0) - f^{\inf}\)}{\varepsilon}
    $$
    seconds, the method guarantees to find an $\epsilon$-stationary point.
\end{lemma}

Now we are ready to prove the theorem.

\begin{restate-theorem}{\ref{theorem_mindflayer_time_complexity}}
    Assume that Assumptions~\ref{ass:lipschitz_constant}, \ref{ass:lower_bound} and \ref{ass:stochastic_variance_bounded} hold. 
    Let 
    $$
        B = \sum_{i=1}^n p_i B_i \quad\text{and}\quad \gamma = \frac{1}{2L}\min\left\{1, \frac{\varepsilon B}{\sigma^2}\right\}
    $$
    in \Cref{alg_mind_flayer}.
    Let $t = \(t_1,\ldots,t_n\)$, $t_1,\ldots,t_n \ge0$.
    Without loss of generality assume that $0< \tau_1 + t_1 \le \cdots \le \tau_n + t_n$.
    Let 
    $$
    t(m) = \(\sum_{j=1}^{m}\frac{p_j}{\tau_j + t_j}\)^{-1}\(S+\sum_{j=1}^{m}p_j\),
    $$
    where $S = \max\left\{ 1, \frac{\sigma^2}{\varepsilon} \right\}$.
    Let $m^* = \arg\min_{m\in [n]} t(m)$, if there are several minimizers we take the smallest one.
    Put 
    $$
    B_i =  \lceil b_i \rceil, 
    \quad     
    b_i = 
    \begin{cases}
    \frac{t(m^*)}{\tau_i + t_i} - 1, & \text{if} \enspace i\le m^*,\\
    0, & \text{if} \enspace i > m^*.
    \end{cases}
    $$
    Then, \algname{MindFlayer SGD} guarantees to find an $\epsilon$-stationary point after
    $$
    T_{\algname{MindFlayer SGD}}(t) 
    \le 8 \times \min_{m\in [n]} \left\{ \(\frac{1}{m}\sum_{j=1}^{m}\frac{p_j}{\tau_j + t_j}\)^{-1}\(\frac{S}{m} + \frac{1}{m}\sum_{j=1}^{m}p_j\)\frac{\Delta L}{\varepsilon} \right\}
    $$
    seconds.
\end{restate-theorem}

\begin{proof}
    First we show that $B_i$-s are valid choice, i.e. $b_i>0$ for $i\le m^*$.
    If $m^* = 1$, then $t(1) = \frac{\tau_1 + t_1}{p_1}(S + p_1)$, thus $b_1 = \frac{S}{p_1} > 0$.
    If $m^* > 1$, then, by its definition, $t(m^*) < t(m^* - 1)$. 
    This implies
    $$
    \left(\sum_{j=1}^{m^*} \frac{p_j}{\tau_j + t_j}\right)^{-1} \(S + \sum_{j=1}^{m^*}p_j\) < 
    \left(\sum_{j=1}^{m^* - 1} \frac{p_j}{\tau_j + t_j}\right)^{-1} \(S + \sum_{j=1}^{m^* - 1}p_j\),
    $$
    leading to
    $$
    \left(\sum_{j=1}^{m^* - 1} \frac{p_j}{\tau_j + t_j}\right) \(S + \sum_{j=1}^{m^*}p_j\) <
    \left(\sum_{j=1}^{m^*} \frac{p_j}{\tau_j + t_j}\right) \(S + \sum_{j=1}^{m^*-1}p_j\)
    $$
    and
    $$
    p_{m^*} \left(\sum_{j=1}^{m^*} \frac{p_j}{\tau_j + t_j}\right) 
    < \frac{p_{m^*}}{\tau_{m^*} + t_{m^*}} \(S + \sum_{j=1}^{m^*}p_j\).
    $$
    From the last inequality, we get that $\tau_{m^*} + t_{m^*} < t(m^*)$, thus $b_i \geq b_{m^*} > 0$ for all $i \leq m^*$. 

    It remains to find the time complexity with these choices of $B_i$.
    From \Cref{lemma_mindflayer_time_complexity}, we have that the time complexity is
    $$
    T_{\algname{MindFlayer SGD}}(t) 
    \le \max_{i \in [n]}\{ B_i\(\tau_i + t_i\) \} \max\left\{ 1, \frac{\sigma^2}{\varepsilon B} \right\} \frac{8 \Delta L}{\varepsilon}.
    $$
    Then,
    $$
    \max_{i \in [n]}\{B_i\(\tau_i + t_i\) \} \le \max_{b_i \ne 0}\{\(b_i + 1\)\(\tau_i + t_i\) \} = t(m^*).
    $$
    On the other hand
    \begin{eqnarray*}
        B 
        &=& \sum_{i=1}^{n} p_i B_i \ge \sum_{i=1}^{n} p_i b_i = \sum_{i=1}^{m^*} t(m^*) \frac{p_i}{\tau_i + t_i} - \sum_{i=1}^{m^*} p_i\\
        &=& \(\sum_{j=1}^{m^*}\frac{p_j}{\tau_j + t_j}\)^{-1}\(S+\sum_{j=1}^{m^*}p_j\)\sum_{i=1}^{m^*} \frac{p_i}{\tau_i + t_i} - \sum_{i=1}^{m^*} p_i = S \ge \frac{\sigma^2}{\varepsilon}.
    \end{eqnarray*}
    Therefore, the time complexity is 
    \begin{eqnarray*}
        T_{\algname{MindFlayer SGD}}(t) 
        &\le& t(m^*) \frac{8 \Delta L}{\varepsilon}\\
        &=& \min_{m\in [n]} \left\{ \(\sum_{j=1}^{m}\frac{p_j}{\tau_j + t_j}\)^{-1}\(S+\sum_{j=1}^{m}p_j\) \right\} \frac{8 \Delta L}{\varepsilon}.
    \end{eqnarray*}
\end{proof}

\subsubsection{Proof of \Cref{lemma_mindflayer_time_complexity}} 
\label{proof_lemma_mindflayer_time_complexity}

\begin{restate-lemma}{\ref{lemma_mindflayer_time_complexity}}
    Assume that Assumptions~\ref{ass:lipschitz_constant}, \ref{ass:lower_bound} and \ref{ass:stochastic_variance_bounded} hold. 
    Let $B = \sum_{i=1}^n p_i B_i$ and $\gamma = \frac{1}{2L}\min\left\{1, \frac{\varepsilon B}{\sigma^2}\right\}$ in Method~\ref{alg_mind_flayer}. 
    Then after 
    $$
    T_{\algname{MindFlayer SGD}}(t) 
    \le \max_{i \in [n]}\{ B_i\(\tau_i + t_i\) \} \max\left\{ 1, \frac{\sigma^2}{\varepsilon B} \right\} \frac{8 L \(f(x_0) - f^{\inf}\)}{\varepsilon}
    $$
    seconds, the method guarantees to find an $\epsilon$-stationary point.
\end{restate-lemma}

\begin{proof}
Let $T_i^j(t_i)$ be the random time taken by client $i$ in the $j$-th attempt of calculating gradient estimator.
We have 
\begin{equation}
    T_i^j(t_i) = 
    \begin{cases}
        \tau_i + \eta_i^j, & \text{if}\enspace \eta_i^j \le t_i, \\
        \tau_i + t_i,      & \text{if}\enspace \eta_i^j > t_i.
    \end{cases}
\end{equation}
Thus, the random time taken for client $i$ to finish it's all $b_i$ trials is 
\begin{equation}
    \cT_i(t_i) := \sum_{j=1}^{B_i} T_i^j(t_i) \le B_i\(\tau_i + t_i\).
\end{equation}
Finally, let $\cT$ be the random time required for one iteration of \algname{MindFlayer SGD}.
We get 
\begin{equation}
    \cT = \max_{i \in [n]} \cT_i(t_i) \le \max_{i \in [n]}\{ B_i\(\tau_i + t_i\) \}.
\end{equation}

It remains to multiply $\cT$ with the number of iterations $K$ given by \Cref{theorem_mindflayer_convergence}.
\end{proof}

\section{Proofs for Heterogeneous Regime}

\subsection{Proof of Theorem~\ref{theorem_heteromindflayer_convergence}}
\label{proof_heteromindflayer_convergence}

Here, we rewrite \algname{Vecna SGD} (\Cref{alg_heteromindflayer}) in a classical SGD way.

\begin{algorithm}[H]
    \caption{\algname{Vecna SGD}}
    \label{alg_heteromindflayer_sgd}

    \begin{algorithmic}[1]
    \STATE \textbf{Input:} starting point $x^0$, stepsize $\gamma$, time budgets $t_1,\ldots,t_n\ge 0$, batch sizes $b_1,\ldots,b_n \ge 0$,
    \FOR{$k = 0, 1, \dots, K - 1$}
    \STATE $g^k = \frac{1}{n}\sum_{i=1}^{n} \frac{1}{p_i B_i}\sum_{j=1}^{B_i}I\(\eta_i^j \le t_i\) \nabla f_i\(x^k;\xi_i^j\)$
    \STATE $x^{k+1} = x^k - \gamma g^{k}$
    \ENDFOR
    \end{algorithmic}
\end{algorithm}
where $p_i = F(t_i) = P(\eta_i \le t_i)$.

To prove the theorem we need to establish some properties of the gradient estimator.
First, we need an unbiased estimator.

\begin{lemma}[Proof in \Cref{proof_lemma_heteromindflayer_unbiased}]
    \label{lemma_heteromindflayer_unbiased}
    The gradient estimator in \Cref{alg_heteromindflayer_sgd} given by 
    $$
    g(x) := \frac{1}{n}\sum_{i=1}^{n} \frac{1}{p_i B_i}\sum_{j=1}^{B_i}I\(\eta_i^j \le t_i\) \nabla f_i\(x;\xi_i^j\)
    $$
    is unbiased, i.e., $\E{g(x)} = \nabla f(x)$ for all $x\in\R^d$.
\end{lemma}

Next, we obtain an upper bound for the variance of this estimator.

\begin{lemma}[Proof in \Cref{proof_lemma_heteromindflayer_variance}]
    \label{lemma_heteromindflayer_variance}
    The gradient estimator in \Cref{alg_heteromindflayer_sgd} given by 
    $$
        g(x) := \frac{1}{n}\sum_{i=1}^{n} \frac{1}{p_i B_i}\sum_{j=1}^{B_i}I\(\eta_i^j \le t_i\) \nabla f_i\(x;\xi_i^j\)
    $$
    satisfies
    $$
        \E{\norm{g(x)^2}} 
        \leq \frac{2L\Delta}{n} \max_{i\in [n]} \left\{ \frac{1-p_i}{p_i B_i} \right\}
            + \sqnorm{\nabla f(x)} 
            + \frac{\sigma^2}{n^2}\sum_{i=1}^{n} \frac{1}{p_i B_i} 
            + \frac{2L}{n} \max_{i\in [n]} \left\{ \frac{1-p_i}{p_i B_i} \right\} \Delta^{\inf}.
    $$
\end{lemma}

We are ready to prove \Cref{theorem_heteromindflayer_convergence}.
First, let us restate the theorem.

\begin{restate-theorem}{\ref{theorem_heteromindflayer_convergence}}
    Assume that Assumptions~\ref{ass:lipschitz_constant}, \ref{ass:lower_bound} and \ref{ass:stochastic_variance_bounded} hold for function $f_i$ for all $i\in [n]$.
    Let 
    $$
    \gamma = \min\left\{\frac{1}{\sqrt{L\alpha K}}, \frac{1}{L \beta}, \frac{\varepsilon}{2 L \zeta}\right\}
    $$
    in \Cref{alg_heteromindflayer}.
    Then after 
    \begin{align*}
        K \geq \frac{12 \Delta L}{\varepsilon}\max\left\{\beta, \frac{12\Delta \alpha}{\varepsilon}, \frac{2 \zeta}{\varepsilon}\right\},
    \end{align*}
    iterations, the method guarantees that $\min_{0\le k \le K} \E{\norm{\nabla f(x^k)}^2} \leq \varepsilon$, and 
    $$
        \alpha = \frac{L}{n} \max_{i\in [n]} \left\{\frac{1-p_i}{p_i B_i}\right\}, \quad 
        \beta = 1, \quad 
        \zeta = \frac{\sigma^2}{n^2}\sum_{i=1}^{n} \frac{1}{p_i B_i} 
            + \frac{2L}{n} \max_{i\in [n]} \left\{ \frac{1-p_i}{p_i B_i} \right\} \Delta^{\inf}.
    $$
\end{restate-theorem}

\begin{proof}
    Note that \Cref{alg_heteromindflayer} can be viewed as a special case of classical stochastic gradient descent (SGD), as reformulated in \Cref{alg_heteromindflayer_sgd}.
    We need to verify that the gradient estimator fulfills the conditions required by classical SGD (\Cref{theorem_sgd}). 
    The two preceding lemmas address this requirement precisely. 
    Specifically, \Cref{lemma_heteromindflayer_unbiased} confirms that the gradient estimator used in \Cref{alg_heteromindflayer_sgd} is unbiased, while \Cref{lemma_heteromindflayer_variance} verifies that the variance of this estimator meets the conditions specified in \Cref{ass:ABC}.
    Consequently, it remains to apply \Cref{theorem_sgd}.
\end{proof}

\subsubsection{Proof of \Cref{lemma_heteromindflayer_unbiased}}
\label{proof_lemma_heteromindflayer_unbiased}

\begin{restate-lemma}{\ref{proof_lemma_heteromindflayer_unbiased}}
    The gradient estimator in \Cref{alg_heteromindflayer_sgd} given by 
    $$
    g(x) := \frac{1}{n}\sum_{i=1}^{n} \frac{1}{p_i B_i}\sum_{j=1}^{B_i}I\(\eta_i^j \le t_i\) \nabla f_i\(x;\xi_i^j\)
    $$
    is unbiased, i.e., $\E{g(x)} = \nabla f(x)$ for all $x\in\R^d$.
\end{restate-lemma}

\begin{proof}
    This follows from direct computation:
    \begin{eqnarray*}
        \E{g(x)} 
        &=& \E{\frac{1}{n}\sum_{i=1}^{n} \frac{1}{p_i B_i}\sum_{j=1}^{B_i}I\(\eta_i^j \le t_i\) \nabla f_i\(x;\xi_i^j\)}\\
        &=& \frac{1}{n}\sum_{i=1}^{n} \frac{1}{p_i B_i}\sum_{j=1}^{B_i}\E{I\(\eta_i^j \le t_i\) \nabla f_i\(x;\xi_i^j\)}\\
        &\overset{\(\eta_i^j \ind \xi_i^j\)}{=}& \frac{1}{n}\sum_{i=1}^{n} \frac{1}{p_i B_i}\sum_{j=1}^{B_i}\E{I\(\eta_i^j \le t_i\)} \E{\nabla f_i\(x;\xi_i^j\)}\\
        &=& \frac{1}{n}\sum_{i=1}^{n} \frac{1}{p_i B_i}\sum_{j=1}^{B_i}p_i \nabla f_i(x)\\
        &=& \frac{1}{n}\sum_{i=1}^{n} \nabla f_i(x) \\
        &=& \nabla f(x).
    \end{eqnarray*}
\end{proof}

\subsubsection{Proof of \Cref{lemma_heteromindflayer_variance}}
\label{proof_lemma_heteromindflayer_variance}

\begin{restate-lemma}{\ref{lemma_heteromindflayer_variance}}
    The gradient estimator in \Cref{alg_heteromindflayer_sgd} given by 
    $$
    g(x) := \frac{1}{n}\sum_{i=1}^{n} \frac{1}{p_i B_i}\sum_{j=1}^{B_i}I\(\eta_i^j \le t_i\) \nabla f_i\(x;\xi_i^j\)
    $$
    satisfies
    $$
        \E{\norm{g(x)^2}} 
        \leq \frac{2L\Delta}{n} \max_{i\in [n]} \left\{ \frac{1-p_i}{p_i B_i} \right\}
            + \sqnorm{\nabla f(x)} 
            + \frac{\sigma^2}{n^2}\sum_{i=1}^{n} \frac{1}{p_i B_i} 
            + \frac{2L}{n} \max_{i\in [n]} \left\{ \frac{1-p_i}{p_i B_i} \right\} \Delta^{\inf}.
    $$
\end{restate-lemma}

\begin{proof}
    Since $\eta_i^j$ and $\xi_i^j$ are independent of each other for all $i\in [n]$ and $j$, we have
\begin{eqnarray*}
\Var{g(x)} &=& \frac{1}{n^2}\sum_{i=1}^{n} \frac{1}{p_i^2 B_i^2}\sum_{j=1}^{B_i}\Var{I\(\eta_i^j \le t_i\) \nabla f_i\(x;\xi_i^j\)},
\end{eqnarray*}
then we use the fact that 
$$
\Var{XY} = \Var{X}\Var{Y} + \Var{X}\E{Y}^2 + \Var{Y}\E{X}^2,
$$
where $X$ and $Y$ are independent random variables.
Hence, we obtain the following bound on the variance 
\begin{eqnarray*}
  \Var{I\(\eta_i^j \le t_i\) \nabla f_i\(x;\xi_i^j\)} 
  &\le& p_i(1-p_i) \sigma^2 + p_i\(1-p_i\) \sqnorm{\nabla f_i (x)} + \sigma^2 p_i^2\\
  &=& p_i \sigma^2 + p_i\(1-p_i\) \sqnorm{\nabla f_i (x)}.
\end{eqnarray*}
As a result, the variance of $g(x)$ is bounded by
\begin{eqnarray*}
\Var{g(x)} 
&\le& \frac{1}{n^2}\sum_{i=1}^{n} \frac{1}{p_i^2 B_i^2}\sum_{j=1}^{B_i}\( p_i \sigma^2 + p_i\(1-p_i\) \sqnorm{\nabla f_i (x)} \) \\
&=& \frac{1}{n^2}\sum_{i=1}^{n} \frac{1}{p_i B_i}\( \sigma^2 + \(1-p_i\)\sqnorm{\nabla f_i (x)} \).
\end{eqnarray*}

Finally

\begin{eqnarray*}
    \E{\norm{g(x)^2}} 
    &=& \Var{g(x)} + \|\E{g(x)}\|^2\\
    &\overset{\le}{}& \sqnorm{\nabla f(x)} + \frac{1}{n^2}\sum_{i=1}^{n} \frac{1-p_i}{p_i B_i} \sqnorm{\nabla f_i (x)} + \frac{\sigma^2}{n^2}\sum_{i=1}^{n} \frac{1}{p_i B_i}.
\end{eqnarray*}

Next we use $\sqnorm{\nabla f_i (x)} \le 2L \(f_i(x) - f_i^{\inf}\)$, thus

\begin{eqnarray*}
    \E{\norm{g(x)^2}} 
        &\le& \frac{2L}{n^2}\sum_{i=1}^{n} \frac{1-p_i}{p_i B_i} \(f_i(x) - f_i^{\inf}\) + \sqnorm{\nabla f(x)} + \frac{\sigma^2}{n^2}\sum_{i=1}^{n} \frac{1}{p_i B_i} \\
        &\le& \frac{2L}{n} \max_{i\in [n]} \left\{ \frac{1-p_i}{p_i B_i} \right\} \frac{1}{n}\sum_{i=1}^{n} \(f_i(x) - f_i^{\inf}\) + \sqnorm{\nabla f(x)} + \frac{\sigma^2}{n^2}\sum_{i=1}^{n} \frac{1}{p_i B_i} \\
        &=& \frac{2L}{n} \max_{i\in [n]} \left\{ \frac{1-p_i}{p_i B_i} \right\} \(f(x) - f^{\inf}\) + \sqnorm{\nabla f(x)} \\
            && + \; \frac{\sigma^2}{n^2}\sum_{i=1}^{n} \frac{1}{p_i B_i}  + \frac{2L}{n} \max_{i\in [n]} \left\{ \frac{1-p_i}{p_i B_i} \right\} \frac{1}{n}\sum_{i=1}^{n} \(f^{\inf} - f_i^{\inf}\).
\end{eqnarray*}

\end{proof}

\subsection{Proof of Theorem~\ref{theorem_heteromindflayer_time_complexity}} 
\label{proof_theorem_heteromindflayer_time_complexity}

The following lemma gives time complexity for any choice of $B_1,\ldots, B_n$ and $t = \(t_1,\ldots,t_n\)$ in \algname{Vecna SGD}.

\begin{lemma}[Proof in \Cref{proof_lemma_vecna_time_complexity}]
    \label{lemma_vecna_time_complexity}
      Assume that Assumptions~\ref{ass:lipschitz_constant}, \ref{ass:lower_bound} and \ref{ass:stochastic_variance_bounded} hold for function $f_i$ for all $i\in [n]$.
    Let $\gamma = \min\left\{\frac{1}{\sqrt{L\alpha K}}, \frac{1}{L}, \frac{\varepsilon}{2 L \zeta}\right\}$ in \Cref{alg_heteromindflayer}.
    Then after 
    $$
    T_{\algname{Vecna SGD}}(t) 
    \le \max_{i \in [n]}\{ B_i\(\tau_i + t_i\) \} \frac{12 \Delta L }{\varepsilon}\max\left\{1, \frac{12\Delta \alpha}{\varepsilon}, \frac{2 \zeta}{\varepsilon}\right\}
    $$
    seconds, where the method guarantees to find an $\epsilon$-stationary point, where $\Delta = f(x_0) - f^{\inf}$ and
    $$
        \alpha = \frac{L}{n} \max_{i\in [n]} \left\{\frac{1-p_i}{p_i B_i}\right\}, \quad 
        \zeta = \frac{\sigma^2}{n^2}\sum_{i=1}^{n} \frac{1}{p_i B_i} 
            + \frac{2L}{n} \max_{i\in [n]} \left\{ \frac{1-p_i}{p_i B_i} \right\} \Delta^{\inf}.
    $$
\end{lemma}

Now we are ready to prove the theorem.

\begin{restate-theorem}{\ref{theorem_heteromindflayer_time_complexity}}
    Assume that Assumptions~\ref{ass:lipschitz_constant}, \ref{ass:lower_bound} and \ref{ass:stochastic_variance_bounded} hold for function $f_i$ for all $i\in [n]$.
  Let $\gamma = \min\left\{\frac{1}{\sqrt{L\alpha K}}, \frac{1}{L}, \frac{\varepsilon}{2 L }\right\}$ in \Cref{alg_heteromindflayer}, where
  $$
    \alpha = \frac{L}{n} \max_{i\in [n]} \left\{\frac{1-p_i}{p_i B_i}\right\}, \quad 
    \zeta = \frac{\sigma^2}{n^2}\sum_{i=1}^{n} \frac{1}{p_i B_i} 
        + \frac{2L}{n} \max_{i\in [n]} \left\{ \frac{1-p_i}{p_i B_i} \right\} \Delta^{\inf}.
  $$
  Let $t = \(t_1,\ldots,t_n\)$, $t_1,\ldots,t_n \ge0$.
  Without loss of generality assume that $0< \tau_1 + t_1 \le \cdots \le \tau_n + t_n$.
  Let 
  $$
    T = \tau_n + t_n + \[\frac{1}{n}\sum_{i=1}^n \frac{\tau_i + t_i}{p_i}\] \frac{\sigma^2}{n \varepsilon} 
        + \max_{i\in [n]} \left\{\frac{1-p_i}{p_i}\(\tau_i + t_i\)\right\} \frac{L\(\Delta + \Delta^{\inf}\)}{n\varepsilon}.
  $$
  Put 
  $$
  B_i =  \lceil b_i \rceil, 
  \quad     
  b_i = \frac{T}{\tau_i + t_i}.
  $$
  Then, \algname{Vecna SGD} guarantees to find an $\epsilon$-stationary point after
  $$
  T_{\algname{Vecna SGD}}(t) 
  \le 288 \times \frac{\Delta L}{\varepsilon} \(\tau_n + t_n + \[\frac{1}{n}\sum_{i=1}^n \frac{\tau_i + t_i}{p_i}\] \frac{\sigma^2}{n \varepsilon} + \[\frac{1}{n}\sum_{i=1}^n \frac{1-p_i}{p_i}\(\tau_i + t_i\)\] \frac{\Delta L}{n \varepsilon}\)
  $$
  seconds. 
\end{restate-theorem}

\begin{proof}
  Since we have $b_i\ge1$ for all $i\in[n]$, we get
  $$
  \max_{i \in [n]}\{ B_i\(\tau_i + t_i\) \} \le \max_{b_i \ne 0}\{\(b_i + 1\)\(\tau_i + t_i\) \} \le 2 \max_{i \in [n]}\{b_i\(\tau_i + t_i\) \} = 2T.
  $$
  It remains to apply \cref{lemma_vecna_time_complexity}.
  We get 
  \begin{eqnarray*}
    \frac{12\Delta \alpha}{\varepsilon} 
    &=& \frac{12\Delta L}{\varepsilon n} \max_{i\in [n]} \left\{\frac{1-p_i}{p_i B_i}\right\}
        \le \frac{12\Delta L}{\varepsilon n} \max_{i\in [n]} \left\{\frac{1-p_i}{p_i b_i}\right\}\\
    &=& \frac{12\Delta L}{n \varepsilon} \frac{1}{T} \max_{i\in [n]} \left\{\frac{1-p_i}{p_i}\(\tau_i + t_i\)\right\}
        \le 12,
\end{eqnarray*}
  and
  \begin{eqnarray*}
    \frac{2\zeta}{\varepsilon} 
    &=& \frac{2\sigma^2}{\varepsilon n^2}\sum_{i=1}^{n} \frac{1}{p_i B_i} 
        + \frac{4L\Delta^{\inf}}{n \varepsilon} \max_{i\in [n]} \left\{\frac{1-p_i}{p_i B_i}\right\} \\
    &\le& \frac{2\sigma^2}{\varepsilon n^2}\sum_{i=1}^{n} \frac{1}{p_i b_i} 
        + \frac{4L\Delta^{\inf}}{n \varepsilon} \max_{i\in [n]} \left\{\frac{1-p_i}{p_i b_i}\right\} \\
    &\le& \frac{2\sigma^2}{n \varepsilon}\frac{1}{T} \frac{1}{n}\sum_{i=1}^{n} \frac{\tau_i + t_i}{p_i} 
        + \frac{4L\Delta^{\inf}}{n \varepsilon}\frac{1}{T} \max_{i\in [n]} \left\{\frac{1-p_i}{p_i}\(\tau_i + t_i\)\right\}\\
    &\le& 4.
  \end{eqnarray*}

  Finally, we get that \cref{alg_heteromindflayer} returns a solution after 
  \begin{eqnarray*}
  T_{\algname{MindFlayer SGD}}(t) 
  &\le& \max_{i \in [n]}\{ B_i\(\tau_i + t_i\) \} \frac{12 \Delta L}{\varepsilon}\max\left\{1, \frac{12\Delta \alpha }{\varepsilon}, \frac{2 \zeta}{\varepsilon}\right\}\\
  &\le& 288 \frac{\Delta L}{\varepsilon}T\\
  &\le& 288 \frac{\Delta L}{\varepsilon} 
    \(\tau_n + t_n + \[\frac{1}{n}\sum_{i=1}^n \frac{\tau_i + t_i}{p_i}\] \frac{\sigma^2}{n \varepsilon} 
    + \max_{i\in [n]} \left\{\frac{1-p_i}{p_i}\(\tau_i + t_i\)\right\} \frac{L\(\Delta + \Delta^{\inf}\)}{n\varepsilon}\)
  \end{eqnarray*}
  seconds.
\end{proof}

\subsubsection{Proof of \Cref{lemma_vecna_time_complexity}} 
\label{proof_lemma_vecna_time_complexity}

\begin{restate-lemma}{\ref{lemma_vecna_time_complexity}}
    Assume that Assumptions~\ref{ass:lipschitz_constant}, \ref{ass:lower_bound} and \ref{ass:stochastic_variance_bounded} hold for function $f_i$ for all $i\in [n]$.
    Let $\gamma = \min\left\{\frac{1}{\sqrt{L\alpha K}}, \frac{1}{L}, \frac{\varepsilon}{2 L \zeta}\right\}$ in \Cref{alg_heteromindflayer}.
    Then after 
    $$
    T_{\algname{Vecna SGD}}(t) 
    \le \max_{i \in [n]}\{ B_i\(\tau_i + t_i\) \} \frac{12 \Delta L }{\varepsilon}\max\left\{1, \frac{12\Delta \alpha}{\varepsilon}, \frac{2 \zeta}{\varepsilon}\right\}
    $$
    seconds, where the method guarantees to find an $\epsilon$-stationary point, and
    $$
    \alpha = \frac{L}{n} \max_{i\in [n]} \left\{\frac{1-p_i}{p_i B_i}\right\}, \quad 
    \zeta = \frac{\sigma^2}{n^2}\sum_{i=1}^{n} \frac{1}{p_i B_i} 
        + \frac{2L}{n} \max_{i\in [n]} \left\{ \frac{1-p_i}{p_i B_i} \right\} \Delta^{\inf}.
    $$
\end{restate-lemma}

\begin{proof}
Let $T_i^j(t_i)$ be the random time taken by client $i$ in the $j$-th attempt of calculating gradient estimator.
We have 
\begin{equation}
    T_i^j(t_i) = 
    \begin{cases}
        \tau_i + \eta_i^j, & \text{if}\enspace \eta_i^j \le t_i, \\
        \tau_i + t_i,      & \text{if}\enspace \eta_i^j > t_i.
    \end{cases}
\end{equation}
Thus, the random time taken for client $i$ to finish it's all $B_i$ trials is 
\begin{equation}
    \cT_i(t_i) := \sum_{j=1}^{B_i} T_i^j(t_i) \le B_i\(\tau_i + t_i\).
\end{equation}
Finally, let $\cT$ be the random time required for one iteration of \algname{Vecna SGD}.
We get 
\begin{equation}
    \cT = \max_{i \in [n]} \cT_i(t_i) \le \max_{i \in [n]}\{ B_i\(\tau_i + t_i\) \}.
\end{equation}

It remains to multiply $\cT$ with the number of iterations $K$ given by \cref{theorem_heteromindflayer_convergence}.
\end{proof}

\end{document}